\documentclass[10pt,leqno]{amsart}
\usepackage{amssymb}

\usepackage{latexsym,amssymb,amsthm,amsmath,amscd}

\theoremstyle{plain} 

\newtheorem{theorem}{Theorem}[section]
\newtheorem{lemma}{Lemma}[section]

\newtheorem{proposition}{Proposition}[section]

\newtheorem{corollary}{Corollary}[section]

\theoremstyle{remark}
\newtheorem{remark}{Remark}[section]

\theoremstyle{definition}

\newtheorem{definition}{Definition}[section]
\newtheorem{example}{Example}[section]
\numberwithin{equation}{section}
\def\<{\left < }
\def\>{\right >}
\def\({\left ( }
\def\){\right )}

\def\e{\eqref}
\def\csch{\,{\rm csch}\,}
\def\sech{\,{\rm sech}\,}

\def\x{\text{ $\frac{\partial}{\partial x}$} }

\def\xj{\text{ ${\partial\over{\partial x_j}}$}}

\def\x1{\text{ ${\partial\over{\partial x_1}}$}}

\begin{document}

\title[Hamiltonian-stationary Lagrangian submanifolds]{Construction of Hamiltonian-stationary Lagrangian submanifolds of constant curvature $\varepsilon$ in complex space  forms $\tilde M^n(4\varepsilon)$}

\author[B.-Y. Chen]{ Bang-Yen Chen}

 

 \maketitle

\begin{abstract}   Lagrangian submanifolds of a Kaehler manifold are called Hamil\-tonian-stationary (or $H$-stationary for short) if it is  a critical point of the area functional restricted to compactly supported Hamiltonian variations. 
In \cite{cdvv}  an effective method to constructing  Lagrangian submanifolds of constant curvature $\varepsilon$ in complex space form $M^n(4\varepsilon)$ was introduced.
In this article we survey recent results on construction of Hamiltonian-stationary Lagrangian submanifolds in complex space forms using this method.
 
\end{abstract}

\section{Introduction.}

Let $\,\tilde M^n(4 \varepsilon)\, $ denote the complex projective $n$-space $\, CP^n(4 \varepsilon)$, the complex Euclidean $n$-space ${\bf C}^n$ or the complex hyperbolic $n$-space $CH^n(4 \varepsilon)$ according to $ \varepsilon>0,\,\varepsilon=0$ or $\varepsilon<0$, respectively.  

The Kaehler 2-form $\omega$ is defined by $\omega(\cdot\,,\cdot)=\<J\cdot,\cdot\>$, where $J$ is the complex structure.  An isometric immersion  $\psi:M^n\to \tilde M^n(4  \varepsilon)$ of an $n$-manifold
$M$ into  $\tilde M^n(4  \varepsilon)$ is called {\it Lagrangian\/} if $\psi^*\omega=0$ on $M$.
A vector field $X$ on $\tilde M^n(4  \varepsilon)$ is called Hamiltonian if ${\mathcal L}_X\omega=f\omega$ for some smooth function $f$ on $\tilde M^n(4  \varepsilon)$, where $\mathcal L$ is the Lie derivative. Thus, there exists a smooth real-valued function $\varphi$ on $\tilde M^n(4 \varepsilon)$ such that $X=J\tilde\nabla \varphi$, where $\tilde\nabla$ is the gradient.
The diffeomorphisms of the flux $\phi_t$ of $X$ transform Lagrangian submanifolds into Lagrangian submanifolds.

A normal vector field $\xi$ to a Lagrangian immersion $\psi:M^n\to \tilde M^n(4 \varepsilon)$ is called Hamiltonian if $\xi=J\nabla f$, where $f$ is a smooth function on $M^n$ and $\nabla f$ is the gradient of $f$ with respect to the induced metric.
If $f\in C^\infty_0(M)$ and $\psi_t:M\to \tilde M^n(4  \varepsilon)$ is a variation of $\psi$ with $\psi_0=\psi$ and variational vector field $\xi$, then the first variation of the volume functional is
$$\frac{d}{dt}_{|_{t=0}}{\rm vol}(M,\psi_t^*g)=-\int_M f\,{\rm div} JH dM,$$
where $H$ is the mean curvature vector of  $\psi$ and $div$ is the divergence on $M^n$. Critical points of this variational functional are called $H$-stationary or Hamiltonian-stationary (cf. \cite{Oh}).    Among others,  $H$-stationary Lagrangian submanifolds in complex space forms have been studied in  \cite{a3}-\cite{cd},  \cite{cg}-\cite{Oh}.

An effective method using twisted products for constructing Lagrangian immersions of a real space form $M^n(\varepsilon)$ into a complex space form $\tilde M^n(4\varepsilon)$ was developed by Chen, Dillen, Verstraelen and Vrancken in \cite{cdvv}. 

One main  result of \cite{cdvv} states that if the twistor form of  a twisted product decomposition ${\mathcal TP}_{f_1\ldots f_\ell}^n( \varepsilon)$ of a simply-connected real space form $M^n(\varepsilon)$ of constant curvature $\varepsilon$ is twisted closed, then it admits a ``unique'' adapted Lagrangian immersion: $$L_{f_1\cdots f_\ell}:{\mathcal TP}_{f_1\ldots f_\ell}^n( \varepsilon)\to \tilde M{}^n(4\varepsilon).$$ Conversely, if $L:  M^n(\varepsilon)\to\tilde M^n(4\varepsilon)$ is a non-totally geodesic Lagrangian  immersion, then $M^n(\varepsilon)$ admits a twisted product decomposition with twisted closed twistor form; moreover, the Lagrangian immersion  is given by the adapted Lagrangian  immersion of the twisted product decomposition. A twisted product decomposition of a real space form is called a warped product decomposition if it is a warped product.

In this article we survey recent results concerning construction of Hamiltonian-stationary Lagrangian submanifolds in complex space forms using this effective method of \cite{cdvv}.
    
  \section{Preliminaries.}

\subsection{Basic notation and formulas} 

Let $L:M\to \tilde M^n(4\tilde  \varepsilon)$ be an isometric immersion of a Riemannian $n$-manifold $M$ into  $M^n(4\tilde  \varepsilon)$. Denote the Riemannian connections of $M$ and  $M^n(4\tilde  \varepsilon)$ by $\nabla$ and $\tilde \nabla$, respectively; and by $D$  the connection on the normal bundle of the submanifold. Let $R$ denote the curvature tensor of $\nabla$.

The formulas of Gauss and Weingarten are \begin{align}  \label{2.1}
&\tilde\nabla_X Y =\nabla_X Y + h(X,Y),\\&\label{2.2}\tilde\nabla_X \xi =- A_\xi X + D_X \xi\end{align}
  for tangent vector fields $X,Y$ and normal  vector field  $\xi$. 
    
If $L:M\to \tilde M^n(4\tilde  \varepsilon)$ is a Lagrangian immersion, then the equations of Gauss and
Codazzi are given respectively by
\begin{align} &\label{2.3}\<R(X,Y)Z,W\> =  \< h( X,W),h(Y,Z)\>  -\< h( X,Z),h(Y,W)\>
\\& \notag \hskip.9in +\varepsilon \{\<X,W\>\<Y,Z\>-\<X,Z\>\<Y,W\>\},\\&\label{2.4}
 (\nabla h)(X,Y,Z) = (\nabla h)(Y,X,Z),\end{align}
where 
$$(\nabla h)(X,Y,Z) = D_X h(Y,Z) - h(\nabla_X Y,Z)- h(Y,\nabla_X Z).$$
For the Lagrangian immersion  we also have (cf. \cite{CO})
\begin{align} \label{2.5} &D_X JY = J \nabla_X Y,\\&\label{2.6} \<h(X,Y),JZ\>=\<h(Y,Z),JX\> =\<h(Z,X),JY\>. \end{align}

 At a given point $p$ on the Lagrangian submanifold $M$, the {\it relative null space} $\mathcal N_p$ at $p$ is the subspace of the tangent space $T_pM$ defined by $${\mathcal N}_p=\{X\in T_pM: h(X,Y)=0 \;\forall Y\in T_pM\}.$$ The dimension of ${\mathcal N}_p$ is called the {\it relative nullity} at $p$.

\subsection{Lagrangian  and Legendrian submanifolds.}

We recall a general method from  \cite{R} for constructing  Lagrangian submanifolds via Hopf's fibration.

{\sc Case (1):} $CP^n(4)$. Let $$S^{2n+1}(1)=\left\{(z_1,\ldots,z_{n+1})\in {\bf C}^{n+1} :\left<z,z\right>=1\right\}$$ be the unit hypersphere in ${\bf C}^{n+1}$ centered at the origin.    On $S^{2n+1}(1)$ we consider the canonical Sasakian structure consisting  of $\phi$ induced from the complex structure $J$ of ${\bf C}^{n+1}$ and the structure vector field $\xi=Jx$ with $x$ being  the position vector.  

An isometric immersion $\psi \colon\;M\to S^{2n+1}(1)$ is called   {\it Legendrian\/}  if $\xi$ is normal to $f_*(TM)$ and $\left<\phi(\psi_*(TM)),\psi_*(TM)\right>=0$, where $\left<\,,\,\right>$ denotes the inner product on ${\bf C}^{n+1}$.   The vectors of $S^{2n+1}(1)$ normal to $\xi$ at a point $z$ define the horizontal subspace $\mathcal H_z$ of the Hopf fibration: $$\pi \colon\;S^{2n+1}(1)\to  CP^n(4).$$

Let $\hat\psi\colon\; M\to CP^n(4)$ be a Lagrangian isometric immersion. Then there is an isometric covering map $\tau\colon\; \hat M \to M$ and a  Legendrian   immersion $\psi\colon\; \hat  M\to S^{2n+1}(1)$ such that $\hat\psi(\tau)=\pi(\psi)$.  Hence every Lagrangian immersion can be lifted locally (or globally if we assume the manifold is simply connected) to a  Legendrian immersion of the same Riemannian manifold. 

Conversely, suppose that  $\psi\colon\; \hat M\to S^{2n+1}(1)$ is a  Legendrian  immersion.  Then
$\hat \psi=\pi(\psi)\colon\;  M\to CP^n(4)$ is  a Lagrangian isometric immersion.   Under this correspondence, the second fundamental forms $h^\psi$ and $h^{\hat \psi}$ of $\psi$ and $\hat \psi$ satisfy $\pi_*h^\psi=h^{\hat \psi}$.  We shall denote $h^\psi$ and $h^{\hat \psi}$ simply by $h$.

\vskip.1in
{\sc Case (2):} $CH^n(-4)$.  Consider  the complex number space  ${\bf C}^{n+1}_1$  with the pseudo Euclidean metric: $g_0=-dz_1d\bar z_1 +\sum_{j=2}^{n+1}dz_jd\bar z_j $. Put 
$$H^{2n+1}_1(-1)=\big\{z=(z_1,z_2,\ldots,z_{n+1}):
\left<z,z\right>=-1\big\},$$
where $\left<\;\,,\;\right>$ is the inner product on ${\bf C}^{n+1}_1$ induced from $g_0$.

Put $$T'_z=\{ z\in {\bf C}^{n+1}: \hbox{Re}\,\left<u,z\right>=\hbox{Re}\, \left<u,iz\right>=0\}$$ and $$H_1^1=\{\lambda\in {\bf C}\colon\; |\lambda|=1\}.$$  Then we have an $H^1_1$-action on $H_1^{2n+1}(-1)$,  $z\mapsto \lambda z$ and at each point $z\in H^{2n+1}_1(-1)$, the vector $iz$ is tangent to the
flow of the action. Since the metric $g_0$ is Hermitian, we have $\hbox{Re} \,g_0(iz,iz)=-1$. The
orbit lies in the negative definite plane spanned by $z$ and $iz$. The quotient space $H^{2n+1}_1/\hskip-.04in \sim$, under the identification from the action, is the complex hyperbolic space $CH^n(-4)$ with  holomorphic sectional curvature $-4$,  with the complex structure $J$ induced from the canonical complex
structure $J$ on $\hbox{\bf C}^{n+1}_1$ via the following pseudo-Riemannian submersion: 
$$\pi\colon\; H^{2n+1}_1(-1)\rightarrow  CH^n(-4).$$

 Just as in Case (1), let $g\colon\; M\to CH^n(-4)$ be a Lagrangian isometric immersion. Then there exists an isometric covering map $\tau\colon\; \hat M \to M$, and a Legendrian isometric immersion $f\colon\; \hat M\to H_1^{2n+1}(-1)$ such that $g(\tau)=\pi(f)$.  Hence every Lagrangian immersion can be lifted locally  to a Legendrian immersion. 

Conversely, let $f:\hat M\to H_1^{2n+1}(-1)$ be a Legendrian immersion.  Then $g=\pi(f)\colon\;  M\to CH^n(-4)$ is again a Lagrangian isometric immersion.  Similarly, under this correspondence, the second fundamental forms $h^f$ and $h^g$ of $f$ and $g$ satisfy $\pi_*h^f=h^g$.  We shall also denote $h^f$ and $h^g$ simply by $h$.
 
Assume that $M$ is a submanifold of $S^{2n+1}(1)$ or $H^{2n+1}_1(-1)$. Denote by
$\tilde\nabla$ and $\nabla$ the Levi-Civita connections of $\hbox{\bf C}^{n+1}$ or  $\hbox{\bf C}^{n+1}_1$ and of $M$, respectively. Let $h$ be the second fundamental form of $M$ in $S^{2n+1}(1)$ or  $H^{2n+1}_1(-1)$. Then we have  
\begin{align}\label{2.7}\tilde \nabla_XY=\nabla_XY +h(X,Y)-\varepsilon\left<X,Y\right> x,\end{align}  where $x$ is the position vector of $M$ in  $\hbox{\bf C}^{n+1}$ or in $\hbox{\bf C}^{n+1}_1$; and
$\varepsilon=1$ or $-1$, according to  the ambient space being $\hbox{\bf C}^{n+1}$ or being $\hbox{\bf C}^{n+1}_1$, respectively.

\section{Warped product decompositions and $H$-stationary.}

Let $(M_j,g_j),j=1,\ldots,m,$ be $m$ Riemannian manifolds, $f_i$ a positive function on $ M_1\times\cdots\times M_m$ and $\pi_i: M_1\times\ldots\times M_m\to M_i$ the $i$-th canonical projection for $i=1,\ldots,m.$ The {\it  twisted  product\/}
$${}_{f_1}M_1\times\cdots\times_{f_m}M_m$$ is the product manifold $M_1\times\ldots\times M_m$ equipped with the twisted product metric $g$ defined by
\begin{align} \label{3.1} &g(X,Y)=f^2_1\cdot g_1({\pi_1}_*X,{\pi_1}_*Y)+\cdots +
f^2_m\cdot g_m({\pi_m}_*X,{\pi_m}_*Y).\end{align}

Let $N^{n-\ell}( \varepsilon)$ be an $(n-\ell)$-dimensional real space form of constant  curvature $ \varepsilon$. For $\ell<n-1$ we consider  the following twisted product:
\begin{align} \label{3.2}{}_{f_1}I_1\times\cdots\times_{f_\ell}I_\ell\times_{1}
N^{n-\ell}( \varepsilon)\end{align} with twisted product metric given by
\begin{align} \label{3.3}g=f^2_1dx_1^2+\cdots+f^2_\ell dx_\ell^2+g_0,\end{align}
where $g_0$ is the canonical metric of $N^{n-\ell}( \varepsilon)$ and $I_1,\ldots,I_\ell$ are open intervals.
When $\ell=n-1$, we shall replace $N^{n-1}(\varepsilon)$ by an open interval. 
If the twisted product  is a real-space-form $M^n( \varepsilon)$, it is called  a {\it twisted product decomposition\/} of  $M^n( \varepsilon)$. We denote such a decomposition by ${\mathcal TP}_{f_1\cdots f_\ell}^n( \varepsilon)$.
 
 Coordinates $x_1,\dots,x_n$ on ${\mathcal TP}_{f_1\cdots f_\ell}^n( \varepsilon)$ are called {\it adapted coordinates} if $\partial/\partial x_j$ is tangent to  $I_j$ for $j=1,\dots, \ell$, the last $n-\ell$ coordinate vectors are tangent to $N^{n-\ell}( \varepsilon)$, and if the metric takes the form \e{3.3}.
 
  The {\it twistor form\/} $\Phi({\mathcal TP})$  on ${\mathcal TP}_{f_1\cdots f_\ell }^n( \varepsilon)$  is defined by
\begin{align} \label{3.4}\Phi({\mathcal TP})=f^2_1dx_1+\cdots+f^2_\ell dx_\ell. \end{align}
 The twistor form is called {\it twisted closed\/} if we have (cf. \cite{cdvv})
\begin{align} \label{3.5} \text{$\sum_{i,j=1}^\ell {{\partial f^2_i}\over{\partial x_j}} $}dx_j\wedge dx_i=0.\end{align}
Obviously, if $\ell=1$, the twisted form $\Phi({\mathcal TP})$ is twisted closed automatically.
 
\begin{theorem} \label{T:1} \cite{cdvv} Let ${\mathcal TP}_{f_1\cdots f_\ell }^n( \varepsilon),\,\ell\in [1,n],$ be a  twisted product decomposition of a simply-connected real-space-form $M^n( \varepsilon)$. If the twistor form $\Phi({\mathcal TP})$ is twisted closed, then, up to rigid motions of $\tilde M^n(4 \varepsilon)$, there is a unique Lagrangian immersion:
\begin{align} \label{3.6} L_{f_1\cdots f_\ell }: {\mathcal TP}_{f_1\cdots f_\ell }^n( \varepsilon) \to \tilde M^n(4 \varepsilon),\end{align}
whose second fundamental form satisfies
\begin{equation}\begin{aligned} &\label{3.7} h\Big({\partial\over{\partial x_j}},{\partial\over{\partial x_j}}\Big)=J\xj,\;\;  j=1,\ldots,\ell ;\;\;  h\Big({\partial\over{\partial x_r}},{\partial\over{\partial x_t}}\Big)=0,\; \hbox{otherwise},\end{aligned}\end{equation} for any adapted coordinate system $\{x_1,\ldots,x_n\}$.

 Conversely, if $\; L: M^n(\varepsilon)\to \tilde M^n(4 \varepsilon)\;$ is a non-totally geodesic Lagrangian immersion of a real-space-form
$M^n(\varepsilon)$ of constant curvature $\varepsilon$ into a complex-space-form $\tilde M^n(4 \varepsilon)$, then $M^n(\varepsilon)$ admits an appropriate twisted product decomposition with twisted closed twistor form and, moreover, the Lagrangian immersion $L$
is given by the corresponding adapted Lagrangian immersion of the twisted product. 
\end{theorem}
 
The $H$-stationary condition for the adapted Lagrangian immersion $ L_{f_1\cdots f_\ell}$ have been  computed by Dong and Han in \cite{DH}.

\begin{proposition} \label{P:3.1}  Let $ L_{f_1\cdots f_\ell}: {\mathcal TP}_{f_1\cdots f_\ell}^n(\varepsilon) \to \tilde M^n(4 \varepsilon)$ be an adapted Lagrangian immersion given in Theorem \ref{T:1}. Then $L_{f_1\cdots f_\ell}$ is $H$-stationary if and only if the twistor  functions $f_1,\ldots,f_\ell$ satisfy
\begin{align} \label{3.8} \text{ $ \sum_{j=1}^\ell  \frac{1}{f_j^4}\frac{\partial f^2_j}{\partial x_j}=\sum_{1\leq i\ne j\leq \ell}\frac{1}{f^2_if^2_j}\frac{\partial f^2_i}{\partial x_j}$}.\end{align}
\end{proposition}
  
  An immediate consequence of this proposition is the following
  
  \begin{corollary} \label{C:3.1} \cite{DH} Any adapted Lagrangian immersion $ L_{f f}: {\mathcal TP}_{ff}^n(\varepsilon) \to \tilde M^n(4 \varepsilon)$ $($with $k=2$ and $f_1=f_2=f)$ is $H$-stationary.
\end{corollary}

\begin{definition}A twisted product decomposition ${\mathcal TP}_{f_1\cdots f_\ell}^n( \varepsilon)$ of a real space form $M^n(\varepsilon)$ is called a {\bf warped product decomposition} if $\ell<n$ and $f_1,\ldots,f_\ell$ are independent of the  adapted coordinates $x_1,\dots,x_\ell$.\end{definition}
 
By applying Theorem \ref{T:1} and Proposition \ref{P:3.1}, we also have the following (cf. \cite{cd,cg}).

\begin{proposition} \label{P:3.2}  Let ${\mathcal TP}_{f_1}^n(\varepsilon)$ be a twisted product decomposition of a real space form $M^n(\varepsilon)$ of constant curvature $\varepsilon$. Then the adapted Lagrangian immersion $L_{f_1}:{\mathcal TP}_{f_1}^n(\varepsilon)\to \tilde M^n(4 \varepsilon)$ is $H$-stationary  if and only if ${\mathcal TP}_{f_1}^n(\varepsilon)$ is a warped product decomposition. 
\end{proposition}

 \begin{proposition}  \label{P:3.3}  If $\, {\mathcal TP}_{f_1\cdots f_\ell}^n( \varepsilon)$ is a warped product decomposition of a simply-connected real space form $M^n(\varepsilon)$, then up to rigid motions,  there exists a unique $H$-stationary Lagrangian immersion \begin{align} \label{3.9} L_{f_1\cdots f_\ell}: {\mathcal TP}_{f_1\cdots f_\ell}^n(c) \to \tilde M^n(4 \varepsilon)\end{align} whose second fundamental form satisfies \e{3.7}.
\end{proposition}

\section{$H$-stationary Lagrangian submanifolds arisen from warped product decompositions.}

It follows from Proposition \ref{P:3.3} that each warped decomposition of a real space space of constant curvature $\varepsilon$ gives rise a Hamiltonian-stationary Lagrangian submanifold in a complex space form of constant holomorphic sectional curvature $4\varepsilon$.

Hamiltonian-stationary Lagrangian submanifolds in complex space forms arisen from warped product decompositions have been completely classified by Chen and Dillen in \cite{cd}. 

\begin{theorem} \label{T:4.1} \cite{cd} There exist two families of non-totally geodesic Hamiltonian-stationary Lagrangian submanifolds in ${\bf C}^n$ arisen from warped product decompositions:

{\rm (a)} Flat  Lagrangian submanifolds defined by
\begin{equation}\begin{aligned} &  L(x_1,\ldots,x_n)=\big(a_1 e^{ix_1},\ldots,a_\ell e^{ix_\ell}, x_{\ell+1},\ldots,x_n\big)\end{aligned}\end{equation} with $a_1,\ldots,a_\ell>0$ and $\ell\in [0,n-1]$.
\smallskip

{\rm (b)} Flat  Lagrangian submanifolds defined by
\begin{equation}\begin{aligned} &   L(x_1,\ldots,x_n)=
\\& 
\Bigg(\tfrac{\sqrt{\sqrt{1+4b_1^2}+1}}{\sqrt{2}(1+4b_1^2)^{1/4}}   e^{\frac{i}{2}(1-\sqrt{1-4b_1^2}\,)x_1}x_{\ell+1}, \ldots,
 \tfrac{\sqrt{\sqrt{1+4b_k^2}+1}}{\sqrt{2}(1+4b_k^2)^{1/4}}  e^{\frac{i}{2}(1-\sqrt{1-4b_k^2}\,)x_k} x_{\ell+k},\\& \hskip.5in 
a_{k+1} e^{ix_{k+1}},\ldots, a_\ell e^{ix_\ell},x_{\ell+k+1},\ldots,x_n
\\& \hskip.04in 
 \tfrac{\sqrt{\sqrt{1+4b_1^2}-1}}{\sqrt{2}(1+4b_1^2)^{1/4}}   e^{\frac{i}{2}(1+\sqrt{1-4b_1^2}\,)x_1}x_{\ell+1}, \ldots,
 \tfrac{\sqrt{\sqrt{1+4b_k^2}-1}}{\sqrt{2}(1+4b_k^2)^{1/4}} e^{\frac{i}{2}(1+\sqrt{1-4b_k^2}\,)x_k}  x_{\ell+k}
 \Bigg) \end{aligned}\end{equation} where $b_1,\ldots,b_k, a_{k+1},\ldots, a_\ell>0$ and $ \ell\in [1,n-1]$.
\end{theorem}

\begin{theorem} \label{T:4.2} There exist two families of non-totally geodesic Hamiltonian-stationary Lagrangian submanifolds of constant curvature one in $CP^n(4)$ arisen from warped product decompositions:

{\rm (a)} Lagrangian submanifolds defined by
\begin{equation}\begin{aligned} &\notag \hskip.0in \hat L=
\(\text{\small$ \frac{2 a_1 e^{\frac{i}{2}x_1}\sin \big(\tfrac{x_1}{2} \sqrt{1+4a_1^2}\big)}{\sqrt{1+4a_1^2}}$}\cos x_{\ell+1} ,\ldots,  \text{\small$ \frac{2 a_\ell e^{\frac{i}{2}x_\ell}\sin  \big(\tfrac{x_\ell}{2} \sqrt{1+4a_\ell^2}\big)}{\sqrt{1+4a_\ell^2}}$}\cos x_{\ell+1} ,\right.
\\&\hskip.2in  e^{\frac{i}{2}x_1}\(\cos \text{\small$\Big(\tfrac{1}{2}\sqrt{1+4a_1^2} x_1\Big)-\frac {i }{\sqrt{1+4a_1^2}} \sin \Big(\tfrac{1}{2}\sqrt{1+4a_1^2} x_1\Big)$}\)\sin\theta_{\ell+1}, \ldots,\\&\hskip.1in  e^{\frac{i}{2}x_\ell} \text{\small$\(\cos\Big(\tfrac{1}{2}\sqrt{1+4a_\ell^2} x_\ell\Big)-\frac {i }{\sqrt{1+4a_\ell^2}} \sin \Big(\tfrac{1}{2}\sqrt{1+4a_\ell^2} x_\ell\Big)\)$}\sin \theta_{2\ell}\prod_{j=1}^{\ell-1} \cos \theta_{\ell+j},
\\&\hskip.3in  \sin\theta_{2\ell+1}\prod_{r=\ell+1}^{2\ell}\cos\theta_{r},\ldots,\sin\theta_{n}\prod_{r=\ell+1}^{n-1}\cos\theta_{r}, \prod_{r=\ell+1}^{n}\cos\theta_{r} \Bigg)
\end{aligned}\end{equation} with $\ell\leq \frac{1}{2}(n+1)$ and $a_1,\ldots,a_\ell>0$.
\smallskip

{\rm (b)}  Lagrangian submanifolds defined by
\begin{equation}\begin{aligned} &\notag \hskip.0in \hat L(x_1,\ldots,x_\ell,\theta_1,\ldots,\theta_{\ell-1})=
\Bigg(\text{\small$ \frac{2 a_1 e^{\frac{i}{2}x_1}}{\sqrt{1+4a_1^2}}\sin \big(\tfrac{1}{2}\sqrt{1+4a_1^2} x_1\big)$}\sin\theta_{\ell+1} ,\ldots, \\&\hskip.4in
\text{\small$ \frac{2 a_{\ell-1} e^{\frac{i}{2}x_{\ell-1}}}{\sqrt{1+4a_{\ell-1}^2}}\sin \big(\tfrac{1}{2} \sqrt{1+4a_{\ell-1}^2}x_{\ell-1}\big)$}\sin\theta_{2\ell-1}\prod_{j=1}^{\ell-2}\cos \theta_{\ell+j} ,
 \\&\hskip.7in
\text{\small$ \frac{2 a_\ell e^{\frac{i}{2}x_{\ell}}}{\sqrt{1+4a_{\ell}^2}}\sin \big(\tfrac{1}{2} \sqrt{1+4a_{\ell}^2}x_{\ell}\big)$}\prod_{j=1}^{\ell-1}\cos \theta_{\ell+j} ,
\end{aligned}\end{equation} 
\begin{equation}\begin{aligned} \notag &\hskip.4in  e^{\frac{i}{2}x_1}\text{\small$ \(\frac {  \sin \big(\tfrac{1}{2}\sqrt{1+4a_1^2} x_1\big)}{\sqrt{1+4a_1^2}} +i\cos \big(\tfrac{1}{2}\sqrt{1+4a_1^2} x_1\big)\)$}\sin\theta_{\ell+1}, \ldots,
\\&\hskip.0in  e^{\frac{i}{2}x_{\ell-1}} \text{\small$\(\frac {  \sin \big(\tfrac{1}{2}\sqrt{1+4a_{\ell-1}^2} x_{\ell-1}\big)}{\sqrt{1+4a_{\ell-1}^2}} +i\cos \(\tfrac{1}{2}\sqrt{1+4a_{\ell-1}^2} x_{\ell-1}\)\)$}\sin \theta_{2\ell-1}\prod_{j=1}^{\ell-2} \cos \theta_{\ell+j},
 \\&\hskip.3in  e^{\frac{i}{2}x_\ell} \text{\small$\(\frac {  \sin \big(\tfrac{1}{2}\sqrt{1+4a_\ell^2} x_\ell\big)}{\sqrt{1+4a_\ell^2}} +i\cos \(\tfrac{1}{2}\sqrt{1+4a_\ell^2} x_\ell\)\)$}\prod_{j=1}^{\ell-1} \cos \theta_{\ell+j},
  \Bigg).
\end{aligned}\end{equation} 
with $n=2\ell-1\geq 3$ and $a_1,\ldots,a_\ell>0$.
\end{theorem}

\begin{theorem} \label{T:4.3} There exist twenty-one families of non-totally geodesic Hamiltonian-stationary Lagrangian submanifolds of constant curvature $-1$ in $CH^n(-4)$ arisen from warped product decompositions:

{\rm (1)}  $n=2$ and Lagrangian submanifolds defined by
\begin{equation}\begin{aligned} & \notag \hskip.2in  \hat L(x_1,\theta_2)=\Bigg(\cosh\theta_2,  \frac{2ae^{\frac{i}{2}x_1}}{\sqrt{1+4a^2}}\sin\Big(\tfrac{1}{2}\sqrt{1+4a^2}x_1\Big)\sinh\theta_2,\hbox{ \hfill}
\\&\hskip.1in e^{\frac{i}{2}x_1}\left\{\cos\Big(\tfrac{\sqrt{1+4a^2}}{2}x_1\Big)-\frac{2a}{\sqrt{1+4a^2}}\sin\Big(\tfrac{1}{2}\sqrt{1+4a^2}x_1\Big)\right\}\sinh\theta_2\Bigg),\; a>0.\end{aligned}\end{equation} 

\smallskip
{\rm (2)}  $n=2$ and Lagrangian submanifolds defined by
\begin{align} & \notag \hskip.2in \hat L=\(e^{\frac{i}{2}x_1} \Big(1-\frac{ix_1}{2}\Big)\cosh \theta_2,\frac{x_1}{2}e^{\frac{i}{2}x_1} \cosh \theta_2,\sinh \theta_2\). \end{align}

\smallskip
{\rm (3)}  $n=2$ and Lagrangian submanifolds defined by
\begin{equation}\begin{aligned}   & \notag\hskip.2in \hat L=\Bigg(e^{\frac{i}{2}x_1}\Big( \cos\big(\tfrac{x_1}{2}\sqrt{1-4a^2}\big)-\frac{i}{\sqrt{1-4a^2}} \sin\big(\tfrac{x_1}{2}\sqrt{1-4a^2}\big)\Big)\cosh \theta_2,
\\&\hskip.7in \frac{2ae^{\frac{i}{2}x_1}}{\sqrt{1-4a^2}}\cosh \theta_2 \sin \(\tfrac{1}{2}\sqrt{1-4a^2}x_1\), \sinh \theta_2\Bigg), \;\; 4a^2<1.\end{aligned}\end{equation}

\smallskip
{\rm (4)}  $n=2$ and Lagrangian submanifolds defined by
\begin{equation}\begin{aligned} & \notag\hskip.2in \hskip.0in \hat L=\Big(e^{\frac{i}{2}x_1}\Big( \cosh\big(\tfrac{x_1}{2}\sqrt{4a^2-1}\big)-\frac{i }{\sqrt{4a^2-1}}\sinh\big(\tfrac{x_1}{2}\sqrt{4a^2-1}\big)\Big)\cosh \theta_2,
\\&\hskip.4in \frac{2ae^{\frac{i}{2}x_1}}{\sqrt{4a^2-1}}\cosh \theta_2 \sinh \(\tfrac{1}{2}\sqrt{4a^2-1}x_1\), \sinh \theta_2 \Big),\;\;  4a^2>1.\end{aligned}\end{equation}

\smallskip
{\rm (5)}  $n=3$ and Lagrangian submanifolds defined by
 \begin{equation}\begin{aligned}  & \notag \hskip.05in  \hat L=
\Bigg( \tfrac{e^{\frac{i}{2}x_1}}{2}(2i+x_1)\cosh \theta_3, \frac{x_1e^{\frac{i}{2}x_1}}{2}\cosh \theta_3,\frac{2b  e^{\frac{i}{2}x_2}\sin\(\frac{x_2}{2}\sqrt{1+4b^2}\,\)}{\sqrt{1+4b^2}}\sinh \theta_3\\& \hskip.3in e^{\frac{i}{2}x_2}\left\{\cos \Big(\tfrac{1}{2}\sqrt{1+4b^2}\,x_2\Big)-  
 \frac{i }{\sqrt{1+4b^2}}\sin\Big(\tfrac{1}{2}\sqrt{1+4b^2}\,x_2\Big)\right\}\sinh \theta_3
\Bigg).
\end{aligned}\end{equation} 

\smallskip
{\rm (6)}  $n=3$ and Lagrangian submanifolds defined by
\begin{equation}\begin{aligned} & \notag \hskip.1in \hat L(x_1,x_2,\theta_3)=\Big(e^{\frac{i}{2}x_1} \cosh\Big(\tfrac{\sqrt{4a^2-1}}{2}x_1\Big)\sin\theta_3,e^{\frac{i}{2}x_1} \sinh\(\tfrac{\sqrt{4a^2-1}}{2}x_1\)\cosh \theta_3, 
\\& \hskip.3in 
e^{\frac{i}{2}x_2} \(\cos \(\tfrac{\sqrt{1+4b^2}}{2}x_2\)-\frac{i }{\sqrt{1+4b^2} } \sin  \(\tfrac{\sqrt{1+4b^2}}{2}x_2\)\) \sinh \theta_3, \\&\hskip.5in  \frac{2be^{\frac{i}{2}x_2} }{\sqrt{1+4b^2}}\sin  \(\tfrac{\sqrt{1+4b^2}}{2}x_2\)\sinh \theta_3\Big),\;\; 4a^2>1.\end{aligned}\end{equation}

\smallskip
{\rm (7)}  $n=3$ and Lagrangian submanifolds defined by
: \begin{equation}\begin{aligned} & \notag\hskip.2in \hat L(x_1,x_2,\theta_3)=
\Bigg(e^{\frac{i}{2}x_1}\(\cos \text{\small$\Big(\tfrac{1}{2}\sqrt{1-4a^2} x_1\Big)-\frac {i  \sin \(\tfrac{1}{2}\sqrt{1-4a^2} x_1\)}{\sqrt{1-4a^2}}$}\)\cosh \theta_3,
\\&\hskip.4in  e^{\frac{i}{2}x_2} \text{\small$\(\cos\(\tfrac{1}{2}\sqrt{1+4b^2} x_2\)-\frac {i  }{\sqrt{1+4b^2}}\sin \(\tfrac{1}{2}\sqrt{1+4b^2} x_2\)\)$}\sinh \theta_3
\\&\hskip.1in \text{\small$ \frac{2 a e^{\frac{i}{2}x_1}\sin \(\frac{1}{2}\sqrt{1-4a^2} x_1\)}{\sqrt{1-4a^2}}$}\cosh \theta_3 , \text{\small$ \frac{2b e^{\frac{i}{2}x_2}\sin \(\frac{1}{2}\sqrt{1+4b^2} x_2\)}{\sqrt{1+4b^2}} $}\sinh \theta_3
  \Bigg), \;\; 4a^2<1.
\end{aligned}\end{equation} 
\smallskip

{\rm (8)} Lagrangian submanifolds defined by
\begin{equation}\begin{aligned} \notag &\hskip.2in \hat L(x_1\ldots,x_{n-1},\theta_n)=
\Bigg(ae^{\theta_n}+\frac{e^{-\theta_n}+2ie^{\theta_n}\sum_{j=1}^{n-1} a_j^2 x_j}{2a},a_1e^{ix_1+\theta_n}, \ldots,\hskip.2in \\&\hskip.2in a_{n-1}e^{ix_{n-1}+\theta_n},\frac{e^{-\theta_n}+2ie^{\theta_n}\sum_{j=1}^{n-1} a_j^2 x_j}{2a} \Bigg)  \end{aligned}\end{equation} 
with $ a=\sqrt{a_1^2+\cdots+a_{n-1}^2}.$

{\rm (9)} $n>2\ell\geq 2$ and Lagrangian submanifolds defined by
\begin{equation}\begin{aligned} & \notag \hskip.2in  \hat L=\sinh\theta_{\ell+1}
\Bigg(\coth \theta_{\ell+1},\text{\small$ \frac{2 a_1 e^{\frac{i}{2}x_1}}{\sqrt{1+4a_1^2}}\sin \big(\tfrac{1}{2}\sqrt{1+4a_1^2} x_1\big)$} \sin \theta_{\ell+2},\ldots,  \\&\hskip.7in
\text{\small$ \frac{2 a_\ell e^{\frac{i}{2}x_{\ell}}}{\sqrt{1+4a_{\ell}^2}}\sin \big(\tfrac{1}{2} \sqrt{1+4a_{\ell}^2}x_{\ell}\big)$}\sin \theta_{2\ell+1}\prod_{j=2}^{\ell}\cos \theta_{\ell+j} ,\\ &\notag \hskip.4in  
e^{\frac{i}{2}x_1}\text{\small$ \( \cos \Big(\tfrac{1}{2}\sqrt{1+4a_1^2} x_1\Big)-\frac {i}{\sqrt{1+4a_1^2}}  \sin \Big(\tfrac{1}{2}\sqrt{1+4a_1^2} x_1\Big) \)$}\sin\theta_{\ell+2}, \ldots, \\& \hskip.0in   e^{\frac{i}{2}x_\ell} \text{\small$ \( \cos \Big(\tfrac{1}{2}\sqrt{1+4a_\ell^2} x_\ell\Big)-\frac {i }{\sqrt{1+4a_\ell^2}}\sin \Big(\tfrac{1}{2}\sqrt{1+4a_\ell^2} x_\ell \Big)\)$}\sin\theta_{2\ell+1}\prod_{j=2}^{\ell} \cos \theta_{\ell+j},
\\&\cos \theta_{\ell+2}\cdots\cos \theta_{2\ell+1}\sin\theta_{2\ell+2},\cdots,\cos \theta_{\ell+2}\cdots\cos \theta_{n-1}\sin\theta_{n} ,
\cos \theta_{\ell+2}\cdots\cos \theta_{n} \Bigg).\end{aligned}\end{equation} 

\smallskip

{\rm (10)} $n>2\ell \geq 2$ and Lagrangian submanifolds defined by

\begin{equation}\begin{aligned} & \notag \hskip.1in  \hat L=\sinh\theta_{\ell+1}\Bigg(e^{\frac{i}{2}x_1}(1-\tfrac{i}{2}x_1)\coth \theta_{\ell+1},\frac{1}{2}e^{\frac{i}{2}x_1}x_1\coth \theta_{\ell+1},\hskip.3in
\\&\hskip.7in\text{\small$ \frac{2 a_2 e^{\frac{i}{2}x_2}}{\sqrt{1+4a_2^2}}\sin \big(\tfrac{1}{2}\sqrt{1+4a_2^2} x_2\big)$} \sin \theta_{\ell+2},\ldots, 
\\&\hskip.7in \text{\small$ \frac{2 a_\ell e^{\frac{i}{2}x_{\ell}}}{\sqrt{1+4a_{\ell}^2}}\sin \big(\tfrac{1}{2} \sqrt{1+4a_{\ell}^2}x_{\ell}\big)$}\sin \theta_{2\ell}\hskip-.05in \prod_{r=\ell+2}^{2\ell-1}\hskip-.05in \cos \theta_{r} ,
\\&\hskip.7in  
e^{\frac{i}{2}x_2}\text{\small$ \( \cos \Big(\tfrac{1}{2}\sqrt{1+4a_2^2} x_2\Big)-\frac {i }{\sqrt{1+4a_2^2}} \sin \Big(\tfrac{1}{2}\sqrt{1+4a_2^2} x_2\Big)\)$}\sin\theta_{\ell+2}, \ldots,
\\&\hskip.3in  e^{\frac{i}{2}x_\ell} \text{\small$ \( \cos \Big(\tfrac{1}{2}\sqrt{1+4a_\ell^2} x_\ell\Big)-\frac {i }{\sqrt{1+4a_\ell^2}} \sin \Big(\tfrac{1}{2}\sqrt{1+4a_\ell^2} x_\ell\Big)\)$}\sin\theta_{2\ell}\prod_{r=\ell+2}^{2\ell-1} \cos \theta_{r},
\\ &
 \hskip.2in \cos \theta_{\ell+2}\cdots\cos \theta_{2\ell}\sin\theta_{2\ell+1},\cdots,
 \cos \theta_{\ell+2}\cdots\cos \theta_{n-1}\sin\theta_{n} ,
\cos \theta_{\ell+2}\cdots\cos \theta_{n}\Bigg).\end{aligned}\end{equation} 
\smallskip

{\rm (11)} $n>2\ell \geq 2$ and Lagrangian submanifolds defined by
\begin{equation}\begin{aligned} &\notag \hskip.1in   \hat L=\sinh\theta_{\ell+1}\Bigg(e^{\frac{i}{2}x_1}\coth \theta_{\ell+1}\Big[\cosh\Big(\tfrac{\sqrt{4a_1^2-1}}{2}x_1\Big)-\tfrac{i \sinh\big(\tfrac{1}{2}\sqrt{4a_1^2-1}x_1\big)}{\sqrt{4a^2-1}}\Big],
\\&\hskip.1in \text{\small$\frac{2a_2e^{\frac{i}{2}x_1}}{\sqrt{4a_1^2-1}}\coth\theta_{\ell+1}\sinh\big(\tfrac{1}{2}\sqrt{4a_1^2-1}x_1\big), 
\frac{2a_2e^{\frac{i}{2}x_2}}{\sqrt{1+4a_2^2}}\sin\Big(\tfrac{1}{2}\sqrt{1+4a_2^2}x_2\Big)\sin\theta_{\ell+2}$},
\\&\hskip.5in e^{\frac{i}{2}x_2}\Big[\text{\small$\cos\Big(\tfrac{1}{2}\sqrt{1+4a_2^2}x_2\Big) - \frac{i }{\sqrt{1+4a_2^2}}  \sin \big(\tfrac{1}{2}\sqrt{1+4a_2^2} x_2\big)$} 
\Big]\sin \theta_{\ell+2},\\& \hskip.5in \ldots,  
\text{\small$ \frac{2 a_\ell e^{\frac{i}{2}x_{\ell}}}{\sqrt{1+4a_{\ell}^2}}\sin \big(\tfrac{1}{2} \sqrt{1+4a_{\ell}^2}x_{\ell}\big)$}\sin\theta_{2\ell}\prod_{r=\ell+2}^{2\ell-1} \cos \theta_{r} ,
 \\&\hskip.1in  e^{\frac{i}{2}x_\ell} \text{\small$ \( \cos \Big(\tfrac{1}{2}\sqrt{1+4a_\ell^2} x_\ell\Big)-\frac {i}{\sqrt{1+4a_\ell^2}}  \sin \Big(\tfrac{1}{2}\sqrt{1+4a_\ell^2} x_\ell\Big)\)$}\sin\theta_{2\ell}
 \hskip-.03in \prod_{r=\ell+2}^{2\ell-1} \cos \theta_{r},\\&\cos \theta_{\ell+2}\cdots\cos \theta_{2\ell}\sin\theta_{2\ell+1},\cdots,\cos \theta_{\ell+2}\cdots\cos \theta_{n-1} \sin\theta_{n},
\cos \theta_{\ell+2}\cdots\cos \theta_{n}\Bigg).\end{aligned}\end{equation}

\smallskip
{\rm (12)} $n>2\ell \geq 2$ and Lagrangian submanifolds defined by
\begin{equation}\begin{aligned} &\notag\hskip.2in \hat L=\sinh\theta_{\ell+1}\Bigg(e^{\frac{i}{2}x_1}\coth \theta_{\ell+1}\Big[\cosh\Big(\tfrac{\sqrt{1-4a_1^2}}{2}x_1\Big)-\tfrac{i \sinh\big(\tfrac{1}{2}\sqrt{1-4a_1^2}x_1\big)}{\sqrt{1-4a^2}}\Big],
\\&\hskip.2in \text{\small$ \frac{2a_1e^{\frac{i}{2}x_1}\coth\theta_{\ell+1}\sinh\big(\tfrac{1}{2}\sqrt{1-4a^2_1}x_1\big)}{\sqrt{1-4a_1^2}}, 
\frac{2a_2e^{\frac{i}{2}x_2}\sin\big(\tfrac{1}{2}\sqrt{1+4a_2^2}x_2\big)}{\sqrt{1+4a_2^2}}\sin\theta_{\ell+2}$},
 \\&\hskip.3in e^{\frac{i}{2}x_2}\text{\small$\Big[\cos\Big(\tfrac{1}{2}\sqrt{1+4a_2^2}x_2\Big) - \frac{i }{\sqrt{1+4a_2^2}} \sin \big(\tfrac{1}{2}\sqrt{1+4a_2^2} x_2\big)
\Big]$}  \sin \theta_{\ell+2}, \ldots, \\& \hskip.5in
\text{\small$ \frac{2 a_\ell e^{\frac{i}{2}x_{\ell}}}{\sqrt{1+4a_{\ell}^2}}\sin \big(\tfrac{1}{2} \sqrt{1+4a_{\ell}^2}x_{\ell}\big)$}\sin\theta_{2\ell}\prod_{r=\ell+2}^{2\ell-1} \cos \theta_{r} ,
\end{aligned}\end{equation} 
\begin{equation}\begin{aligned} \notag\\ &\hskip.1in  e^{\frac{i}{2}x_\ell} \text{\small$ \( \cos \Big(\tfrac{1}{2}\sqrt{1+4a_\ell^2} x_\ell\Big)-\frac {i}{\sqrt{1+4a_\ell^2}} \sin \Big(\tfrac{1}{2}\sqrt{1+4a_\ell^2} x_\ell\Big) \)$}\sin\theta_{2\ell}\prod_{r=\ell+2}^{2\ell-1} \cos \theta_{r}, \\&\cos \theta_{\ell+2}\cdots\cos \theta_{2\ell}\sin\theta_{2\ell+1},\cdots,\cos \theta_{\ell+2}\cdots\cos \theta_{n-1}\sin\theta_{n},
\cos \theta_{\ell+2}\cdots\cos \theta_{n}\Bigg).\end{aligned}\end{equation}

\smallskip
{\rm (13)} $n\geq 3, n>2\ell\geq 6$ and Lagrangian submanifolds defined by
\begin{equation}\begin{aligned} &\notag \hskip.2in   \hat L=\sinh\theta_{\ell+1}
\Bigg(\(1+\tfrac{a_1^2}{2}+ia_1^2x_1\)\coth \theta_{\ell+1}-a_1^2\(\tfrac{1}{2}+ix_1\)\sin\theta_{\ell+2},
\\ &\hskip.1in a_1^2\(\tfrac{1}{2}+ix_1\)\coth\theta_{\ell+1}+
\(1-\tfrac{a_1^2}{2}-ia_1^2x_1\)\sin\theta_{\ell+2},
ia_1 e^{ix_1}(\coth \theta_{\ell+1}-\sin\theta_{\ell+2}),
\\  &\hskip.4in
\text{\small$ \frac{2 a_2 e^{\frac{i}{2}x_{2}}}{\sqrt{1+4a_{2}^2}}\sin \big(\tfrac{1}{2} \sqrt{1+4a_{2}^2}x_{2}\big)$}\cos \theta_{\ell+2}\sin \theta_{\ell+3} ,\\ &\hskip.2in  
e^{\frac{i}{2}x_2}\text{\small$ \( \cos \Big(\tfrac{1}{2}\sqrt{1+4a_2^2} x_2\Big)-\frac {i }{\sqrt{1+4a_2^2}} \sin \Big(\tfrac{1}{2}\sqrt{1+4a_2^2} x_2\Big)\)$}\cos\theta_{\ell+2}\sin\theta_{\ell+3}, 
\\& \hskip1.2in \ldots
\\&\hskip.4in
\text{\small$ \frac{2 a_{\ell} e^{\frac{i}{2}x_{\ell}}}{\sqrt{1+4a_{\ell}^2}}\sin \big(\tfrac{1}{2} \sqrt{1+4a_{\ell}^2}x_{\ell}\big)$}\sin \theta_{2\ell+1}\prod_{j=2}^{\ell}\cos \theta_{\ell+j} ,
\\&\hskip.1in   e^{\frac{i}{2}x_{\ell}} \text{\small$ \( \cos \Big(\tfrac{1}{2}\sqrt{1+4a_{\ell}^2} x_{\ell}\Big)-\frac {i }{\sqrt{1+4a_{\ell}^2}}\sin \Big(\tfrac{1}{2}\sqrt{1+4a_{\ell}^2} x_{\ell}\Big) \)$}\sin\theta_{2\ell+1}\prod_{t=\ell+2}^{2\ell} \cos \theta_{t},
\\&\hskip.0in  \cos\theta_{\ell+2}\cdots\cos\theta_{2\ell+1}\sin\theta_{2\ell+2},\cdots, \cos\theta_{\ell+2}\cdots\cos\theta_{n-1}\sin\theta_{n},
  \cos\theta_{\ell+2}\cdots\cos\theta_n \Bigg).
\end{aligned}\end{equation} 

\smallskip
{\rm (14)} $n=2\ell\geq 4$ and Lagrangian submanifolds defined by
\begin{equation}\begin{aligned} &\notag \hskip.2in   \hat L=\sinh\theta_{\ell+1}
\Bigg(\coth \theta_{\ell+1},\text{\small$ \frac{2 a_1 e^{\frac{i}{2}x_1}}{\sqrt{1+4a_1^2}}\sin \big(\tfrac{1}{2}\sqrt{1+4a_1^2} x_1\big)$} \sin \theta_{\ell+2},
\\ &\hskip.4in  
e^{\frac{i}{2}x_1}\text{\small$ \( \cos \Big(\tfrac{1}{2}\sqrt{1+4a_1^2} x_1\Big)-\frac {i }{\sqrt{1+4a_1^2}} \sin \Big(\tfrac{1}{2}\sqrt{1+4a_1^2} x_1\Big)\)$}\sin\theta_{\ell+2}, 
\\& \hskip1.2in \ldots
 \\&\hskip.4in
\text{\small$ \frac{2 a_{\ell-1} e^{\frac{i}{2}x_{\ell-1}}}{\sqrt{1+4a_{\ell-1}^2}}\sin \big(\tfrac{1}{2} \sqrt{1+4a_{\ell-1}^2}x_{\ell-1}\big)$}\sin \theta_{2\ell}\prod_{j=2}^{\ell-1}\cos \theta_{\ell+j} ,
\\&\hskip.0in   e^{\frac{i}{2}x_{\ell-1}} \text{\small$ \( \cos \(\tfrac{1}{2}\sqrt{1+4a_{\ell-1}^2} x_{\ell-1}\)-\frac {i \sin \(\tfrac{1}{2}\sqrt{1+4a_{\ell-1}^2} x_{\ell-1}\)}{\sqrt{1+4a_{\ell-1}^2}} \)$}\sin\theta_{2\ell}\prod_{t=\ell+2}^{2\ell-1} \cos \theta_{t},
\\ &\hskip.4in
\text{\small$ \frac{2 a_\ell e^{\frac{i}{2}x_{\ell}}}{\sqrt{1+4a_{\ell}^2}}\sin \big(\tfrac{1}{2} \sqrt{1+4a_{\ell}^2}x_{\ell}\big)$}\cos \theta_{\ell+2}\cdots\cos \theta_{2\ell} , \\&\hskip.0in   e^{\frac{i}{2}x_\ell} \text{\small$ \( \cos \Big(\tfrac{1}{2}\sqrt{1+4a_\ell^2} x_\ell\Big)-\frac {i }{\sqrt{1+4a_\ell^2}}\sin \Big(\tfrac{1}{2}\sqrt{1+4a_\ell^2} x_\ell\Big) \)$}\cos \theta_{\ell+2}\cdots\cos \theta_{2\ell}
 \Bigg).\end{aligned}\end{equation} 

\smallskip
{\rm (15)} $n=2\ell\geq 4$ and Lagrangian submanifolds defined by
\begin{equation}\begin{aligned} &\notag\hskip.1in \hat L= \sinh\theta_{\ell+1}\Bigg(\hskip-.02in  e^{\frac{i}{2}x_1}(1\hskip-.02in -\hskip-.02in \tfrac{i}{2}x_1)\coth \theta_{\ell+1},\hskip-.02in  \frac{e^{\frac{i}{2}x_1}x_1}{2}\coth \theta_{\ell+1},
\\&\hskip.0in\text{\small$ \frac{2 a_2 e^{\frac{i}{2}x_2}\sin \big(\tfrac{1}{2}\sqrt{1+4a_2^2} x_2\big)}{\sqrt{1+4a_2^2}}$} \sin \theta_{\ell+2},\ldots,  
\text{\small$ \frac{2 a_\ell e^{\frac{i}{2}x_{\ell}}\sin \big(\tfrac{1}{2} \sqrt{1+4a_{\ell}^2}x_{\ell}\big)}{\sqrt{1+4a_{\ell}^2}}$}\sin \theta_{2\ell}\hskip-.02in  \prod_{r=\ell+2}^{n-1}\hskip-.04in \cos \theta_{r} ,
\\&\hskip.2in  
e^{\frac{i}{2}x_1}\text{\small$ \( \cos \Big(\tfrac{1}{2}\sqrt{1+4a_1^2} x_1\Big)-\frac {i}{\sqrt{1+4a_1^2}} \) \sin \Big(\tfrac{1}{2}\sqrt{1+4a_1^2} x_1\Big)$}\sin\theta_{\ell+2}, \ldots,
\\&\hskip.1in   e^{\frac{i}{2}x_\ell} \Big(\text{\small$ \cos \Big(\tfrac{1}{2}\sqrt{1+4a_\ell^2} x_\ell\Big)-\frac {i }{\sqrt{1+4a_\ell^2}} \sin \Big(\tfrac{1}{2}\sqrt{1+4a_\ell^2} x_\ell\Big)$}\Big) \sin\theta_{n}\prod_{r=\ell+2}^{n-1} \cos \theta_{r}\Bigg). \end{aligned}\end{equation}

\smallskip
{\rm (16)} $n=2\ell\geq 4$ and Lagrangian submanifolds defined by
\begin{equation}\begin{aligned} &\notag\hskip.1in  \hat L= \sinh\theta_{\ell+1}\Bigg(\hskip-.02in  e^{\frac{i}{2}x_1}\coth \theta_{\ell+1}\Big[\cosh\Big(\tfrac{\sqrt{4a_1^2-1}}{2}x_1\Big)\hskip-.02in - \hskip-.02in \tfrac{i \sinh\big(\tfrac{1}{2}\sqrt{4a_1^2-1}x_1\big)}{\sqrt{4a^2-1}}\Big],
\\&\hskip.3in \text{\small$\frac{2a_2e^{\frac{i}{2}x_1}\coth\theta_{\ell+1}}{\sqrt{4a_1^2-1}}\sinh\Big(\tfrac{1}{2}\sqrt{4a_1^2-1}x_1\Big), 
\frac{2a_2e^{\frac{i}{2}x_2}}{\sqrt{1+4a_2^2}}\sin\Big(\tfrac{1}{2}\sqrt{1+4a_2^2}x_2\Big)\sin\theta_{\ell+2},$}
\\&\hskip.5in e^{\frac{i}{2}x_2}\Big[\text{\small$\cos\Big(\tfrac{1}{2}\sqrt{1+4a_2^2}x_2\Big) - \frac{i  }{\sqrt{1+4a_2^2}}\sin \big(\tfrac{1}{2}\sqrt{1+4a_2^2} x_2\big)$} 
\Big]\sin \theta_{\ell+2}, \ldots,  \\&\hskip.7in
\text{\small$ \frac{2 a_\ell e^{\frac{i}{2}x_{\ell}}}{\sqrt{1+4a_{\ell}^2}}\sin \big(\tfrac{1}{2} \sqrt{1+4a_{\ell}^2}x_{\ell}\big)$}\cos \theta_{\ell+2}\cdots\cos\theta_{n-1} \sin\theta_{n} ,
 \\&\hskip.2in  e^{\frac{i}{2}x_\ell} \Big[\text{\small$  \cos \Big(\tfrac{1}{2}\sqrt{1+4a_\ell^2} x_\ell\Big)-\frac {i }{\sqrt{1+4a_\ell^2}} \sin \Big(\tfrac{1}{2}\sqrt{1+4a_\ell^2} x_\ell\Big)$}\Big]\sin\theta_{n}\prod_{r=\ell+2}^{n-1} \cos \theta_{r}\Bigg). \end{aligned}\end{equation}

\smallskip
{\rm (17)} $n=2\ell\geq 4$ and  Lagrangian submanifolds defined by
\begin{equation}\begin{aligned} &\notag\hskip.1in  \hat L=
\sinh\theta_{\ell+1}\Bigg(\hskip-.02in e^{\frac{i}{2}x_1}\coth \theta_{\ell+1}\Big[\cosh\Big(\tfrac{\sqrt{1-4a_1^2}}{2}x_1\Big)\hskip-.02in - \hskip-.02in \tfrac{i \sinh\big(\tfrac{1}{2}\sqrt{1-4a_1^2}x_1\big)}{\sqrt{1-4a^2}}\Big],
\\ &\hskip.1in \text{\small$\frac{2a_1e^{\frac{i}{2}x_1}}{\sqrt{1-4a_1^2}}\coth\theta_{\ell+1}\sinh\Big(\tfrac{x_1}{2}\sqrt{1-4a^2_1}\Big), 
\frac{2a_2e^{\frac{i}{2}x_2}}{\sqrt{1+4a_2^2}}\sin\Big(\tfrac{x_2}{2}\sqrt{1+4a_2^2}\Big)\sin\theta_{\ell+2},$}
\\&\hskip.4in e^{\frac{i}{2}x_2}\Big[\text{\small$\cos\Big(\tfrac{x_2}{2}\sqrt{1+4a_2^2}\Big)$} -\text{\small$ \frac{i  }{\sqrt{1+4a_2^2}}\sin \big(\tfrac{x_2}{2}\sqrt{1+4a_2^2} \big)$} 
\Big]\sin \theta_{\ell+2}, \ldots,  \\&\hskip.5in
\text{\small$ \frac{2 a_\ell e^{\frac{i}{2}x_{\ell}}}{\sqrt{1+4a_{\ell}^2}}\sin \big(\tfrac{1}{2} \sqrt{1+4a_{\ell}^2}x_{\ell}\big)$}\sin\theta_{n} \cos \theta_{\ell+2}\cdots\cos\theta_{n-1} ,
 \\&\hskip.2in  e^{\frac{i}{2}x_\ell} \text{\small$ \Big[\cos \Big(\tfrac{1}{2}\sqrt{1+4a_\ell^2} x_\ell\Big)-\frac {i }{\sqrt{1+4a_\ell^2}}\sin \Big(\tfrac{1}{2}\sqrt{1+4a_\ell^2} x_\ell\Big)$}\Big]\sin\theta_{n}\prod_{r=\ell+2}^{n-1} \cos \theta_{r}\Bigg).\end{aligned}\end{equation}

\smallskip
{\rm (18)} $n=2\ell\geq 6$ and Lagrangian submanifolds defined by
\begin{equation}\begin{aligned} &\notag \hskip.2in   \hat L=\sinh\theta_{\ell+1}
\Bigg(\big(1+\tfrac{a_1^2}{2}+ia_1^2x_1\big)\coth \theta_{\ell+1}-a_1^2\(\tfrac{1}{2}+ix_1\)\sin\theta_{\ell+2},
\end{aligned}\end{equation}
 \begin{equation}\begin{aligned} \notag \\ &\hskip.1in a_1^2\(\tfrac{1}{2}+ix_1\)\coth\theta_{\ell+1}+
\big(1-\tfrac{a_1^2}{2}-ia_1^2x_1\big)\sin\theta_{\ell+2},
ia_1 e^{ix_1}(\coth \theta_{\ell+1}-\sin\theta_{\ell+2}),
\\\ &\hskip.6in
\text{\small$ \frac{2 a_2 e^{\frac{i}{2}x_{2}}}{\sqrt{1+4a_{2}^2}}\sin \big(\tfrac{1}{2} \sqrt{1+4a_{2}^2}x_{2}\big)$}\cos \theta_{\ell+2}\sin \theta_{\ell+3} ,\\ &\hskip.4in  
e^{\frac{i}{2}x_2}\text{\small$ \( \cos \Big(\tfrac{1}{2}\sqrt{1+4a_2^2} x_2\Big)-\frac {i }{\sqrt{1+4a_2^2}} \sin \Big(\tfrac{1}{2}\sqrt{1+4a_2^2} x_2\Big)\)$}\cos\theta_{\ell+2}\sin\theta_{\ell+3}, 
\\& \hskip1.2in \ldots
\\ &\hskip.4in
\text{\small$ \frac{2 a_{\ell-1} e^{\frac{i}{2}x_{\ell-1}}}{\sqrt{1+4a_{\ell-1}^2}}\sin \Big(\frac{x_{\ell-1}}{2} \sqrt{1+4a_{\ell-1}^2}\,\Big)$}\sin \theta_{2\ell}\prod_{j=2}^{\ell-1}\cos \theta_{\ell+j} ,
\\&\hskip.1in   e^{\frac{i}{2}x_{\ell-1}} \text{\small$ \( \cos \(\frac{x_{\ell-1}}{2}\sqrt{1+4a_{\ell-1}^2}\, \)-\frac {i \sin \(\frac{x_{\ell-1}}{2}\sqrt{1+4a_{\ell-1}^2}\, \)}{\sqrt{1+4a_{\ell-1}^2}} \)$}\sin\theta_{2\ell}\prod_{t=\ell+2}^{2\ell-1} \cos \theta_{t},
\\ &\hskip.4in
\text{\small$ \frac{2 a_\ell e^{\frac{i}{2}x_{\ell}}}{\sqrt{1+4a_{\ell}^2}}\sin \big(\tfrac{1}{2} \sqrt{1+4a_{\ell}^2}x_{\ell}\big)$}\cos \theta_{\ell+2}\cdots\cos \theta_{2\ell} , \\&\hskip.2in   e^{\frac{i}{2}x_\ell} \text{\small$ \( \cos \Big(\tfrac{1}{2}\sqrt{1+4a_\ell^2} x_\ell\Big)-\frac {i }{\sqrt{1+4a_\ell^2}}\sin \Big(\tfrac{1}{2}\sqrt{1+4a_\ell^2} x_\ell\Big) \)$}\cos \theta_{\ell+2}\cdots\cos \theta_{2\ell}
 \Bigg).\end{aligned}\end{equation} 

\smallskip
{\rm (19)} $n=2\ell-1\geq 5$ and Lagrangian submanifolds defined by
\begin{equation}\begin{aligned} &\notag\hskip.2in \hat L=\sinh\theta_{\ell+1}\Bigg(e^{\frac{i}{2}x_1}\Big(1-\frac{i}{2}x_1\Big)\coth \theta_{\ell+1},\frac{x_1}{2}e^{\frac{i}{2}x_1}\coth \theta_{\ell+1},
\\&\hskip.5in\text{\small$ \frac{2 a_2 e^{\frac{i}{2}x_2}}{\sqrt{1+4a_2^2}}$}\sin \big(\tfrac{1}{2}\sqrt{1+4a_2^2}\, x_2\big) \sin \theta_{\ell+2},\ldots,  \\&\hskip.7in
\text{\small$ \frac{2 a_{\ell-1} e^{\frac{i}{2}x_{\ell-1}}}{\sqrt{1+4a_{\ell-1}^2}}\sin \big(\tfrac{1}{2} \sqrt{1+4a_{\ell-1}^2}x_{\ell-1}\big)$}\sin \theta_{2\ell-1}\prod_{r=\ell+2}^{2\ell-2}\cos \theta_{r} ,
\\ &\hskip.4in  
e^{\frac{i}{2}x_2}\text{\small$ \( \cos \Big(\tfrac{1}{2}\sqrt{1+4a_2^2} x_2\Big)-\frac {i }{\sqrt{1+4a_2^2}} \sin \(\tfrac{1}{2}\sqrt{1+4a_2^2} x_2\)\)$}\sin\theta_{\ell+2}, \ldots,
 \\&\hskip.3in  e^{\frac{i}{2}x_{\ell-1}} \text{\small$ \( \cos \Big(\tfrac{1}{2}\sqrt{1+4a_{\ell-1}^2} x_{\ell-1}\Big)-\tfrac {i \sin \Big(\tfrac{1}{2}\sqrt{1+4a_{\ell-1}^2} x_{\ell-1}\Big)}{\sqrt{1+4a_{\ell-1}^2}} \)$}\sin\theta_{2\ell-1}\prod_{r=\ell+2}^{2\ell-2} \cos \theta_{r},\\&\hskip.7in
\text{\small$
 \frac{2 a_\ell e^{\frac{i}{2}x_{\ell}}}{\sqrt{1+4a_{\ell}^2}}
 \sin \big(\tfrac{1}{2} \sqrt{1+4a_{\ell}^2}x_{\ell}\big) $}\prod_{r=\ell+2}^{2\ell-1}\cos \theta_{r} ,
 \\&\hskip.3in  e^{\frac{i}{2}x_\ell} \text{\small$ \( \cos \(\tfrac{1}{2}\sqrt{1+4a_\ell^2} x_\ell\)-\frac {i}{\sqrt{1+4a_\ell^2}} \sin \Big(\tfrac{1}{2}\sqrt{1+4a_\ell^2} x_\ell\Big) \)$}\prod_{r=\ell+2}^{2\ell-1} \cos \theta_{r}\Bigg).\end{aligned}\end{equation}

\smallskip
{\rm (20)} $n=2\ell-1\geq 5$ and Lagrangian submanifolds defined by
\begin{equation}\begin{aligned} &\notag\hskip.2in \hat L=\sinh\theta_{\ell+1}\Bigg(e^{\frac{i}{2}x_1}\coth \theta_{\ell+1}\Big[\cosh\Big(\tfrac{\sqrt{4a_1^2-1}}{2}x_1\Big)-\tfrac{i \sinh\big(\tfrac{1}{2}\sqrt{4a_1^2-1}x_1\big)}{\sqrt{4a^2-1}}\Big],
\\&\hskip.1in \text{\small$ \frac{2a_1e^{\frac{i}{2}x_1}}{\sqrt{4a_1^2-1}}\coth\theta_{\ell+1}\sinh\big(\tfrac{1}{2}\sqrt{4a_1^2-1}x_1\big)$}, 
\text{\small$ \frac{2 a_2 e^{\frac{i}{2}x_2}}{\sqrt{1+4a_2^2}}\sin \big(\tfrac{1}{2}\sqrt{1+4a_2^2}\, x_2\big)$} \sin \theta_{\ell+2}, 
\end{aligned}\end{equation}
 \begin{equation}\begin{aligned} \notag \\&\hskip.7in \ldots, 
\text{\small$ \frac{2 a_{\ell-1} e^{\frac{i}{2}x_{\ell-1}}}{\sqrt{1+4a_{\ell-1}^2}}\sin \big(\tfrac{1}{2} \sqrt{1+4a_{\ell-1}^2}x_{\ell-1}\big)$}\sin \theta_{2\ell-1}\prod_{r=\ell+2}^{2\ell-2}\cos \theta_{r} ,
 \\&  \hskip.4in  
e^{\frac{i}{2}x_2}\text{\small$ \Big( \cos \Big(\tfrac{1}{2}\sqrt{1+4a_2^2} x_2\Big)-\frac {i }{\sqrt{1+4a_2^2}} \sin \Big(\tfrac{1}{2}\sqrt{1+4a_2^2} x_2\Big)\Big)$}\sin\theta_{\ell+2}, \ldots,
 \\&\hskip.3in  e^{\frac{i}{2}x_{\ell-1}} \text{\small$ \( \cos \Big(\tfrac{1}{2}\sqrt{1+4a_{\ell-1}^2} x_{\ell-1}\Big)-\tfrac {i \sin \Big(\tfrac{1}{2}\sqrt{1+4a_{\ell-1}^2} x_{\ell-1}\Big)}{\sqrt{1+4a_{\ell-1}^2}} \)$}\sin\theta_{2\ell-1}\prod_{r=\ell+2}^{2\ell-2} \cos\theta_r,
\\&\hskip.7in
\text{\small$ \frac{2 a_\ell e^{\frac{i}{2}x_{\ell}}}{\sqrt{1+4a_{\ell}^2}}\sin \big(\tfrac{1}{2} \sqrt{1+4a_{\ell}^2}x_{\ell}\big)$}\cos \theta_{\ell+2}\cdots\cos\theta_{2\ell-1},
\\ &\hskip.3in  e^{\frac{i}{2}x_\ell} \text{\small$ \( \cos \Big(\tfrac{1}{2}\sqrt{1+4a_\ell^2} x_\ell\Big)-\frac {i}{\sqrt{1+4a_\ell^2}} \sin \Big(\tfrac{1}{2}\sqrt{1+4a_\ell^2} x_\ell\Big) \)$}\prod_{r=\ell+2}^{2\ell-1} \cos \theta_{r}\Bigg).\end{aligned}\end{equation}

\smallskip
{\rm (21)} $n=2\ell-1\geq 5$ and Lagrangian submanifolds defined by
\begin{equation}\begin{aligned} &\notag\hskip.0in \hat L=\sinh\theta_{\ell+1}\Bigg(e^{\frac{i}{2}x_1}\coth \theta_{\ell+1}\Big[\cos\Big(\tfrac{\sqrt{1-4a_1^2}}{2}x_1\Big)-\frac{i }{\sqrt{1-4a^2}}\sin\big(\tfrac{1}{2}\sqrt{1-4a_1^2}x_1\big)\Big],
\\&\hskip.1in \text{\small$ \frac{2a_1e^{\frac{i}{2}x_1}}{\sqrt{1-4a_1^2}}\coth\theta_{\ell+1}\sin\big(\tfrac{1}{2}\sqrt{1-4a_1^2}x_1\big)$}, 
\text{\small$ \frac{2 a_2 e^{\frac{i}{2}x_2}}{\sqrt{1+4a_2^2}}\sin \big(\tfrac{1}{2}\sqrt{1+4a_2^2}\, x_2\big)$} \sin \theta_{\ell+2}, \\&\hskip.7in \ldots, 
\text{\small$ \frac{2 a_{\ell-1} e^{\frac{i}{2}x_{\ell-1}}}{\sqrt{1+4a_{\ell-1}^2}}\sin \big(\tfrac{1}{2} \sqrt{1+4a_{\ell-1}^2}x_{\ell-1}\big)$}\sin \theta_{2\ell-1}\prod_{r=\ell+2}^{2\ell-2}\cos \theta_{r} ,
\\&\hskip.4in  
e^{\frac{i}{2}x_2}\text{\small$ \( \cos \Big(\tfrac{1}{2}\sqrt{1+4a_2^2} x_2\Big)-\frac {i }{\sqrt{1+4a_2^2}} \sin \(\tfrac{1}{2}\sqrt{1+4a_2^2} x_2\)\)$}\sin\theta_{\ell+2}, \ldots,
 \\&\hskip.1in  e^{\frac{i}{2}x_{\ell-1}} \text{\small$ \( \cos \Big(\tfrac{1}{2}\sqrt{1+4a_{\ell-1}^2} x_{\ell-1}\Big)-\tfrac {i \sin \Big(\tfrac{1}{2}\sqrt{1+4a_{\ell-1}^2} x_{\ell-1}\Big)}{\sqrt{1+4a_{\ell-1}^2}} \)$}\sin\theta_{2\ell-1}\prod_{r=\ell+2}^{2\ell-2} \cos \theta_{r},\\&\hskip.7in
\text{\small$ \frac{2 a_\ell e^{\frac{i}{2}x_{\ell}}}{\sqrt{1+4a_{\ell}^2}}\sin \big(\tfrac{1}{2} \sqrt{1+4a_{\ell}^2}x_{\ell}\big)$}\prod_{r=\ell+2}^{2\ell-1}\cos \theta_{r} ,
 \\&\hskip.3in  e^{\frac{i}{2}x_\ell} \text{\small$ \( \cos \Big(\tfrac{1}{2}\sqrt{1+4a_\ell^2} x_\ell\Big)-\frac {i}{\sqrt{1+4a_\ell^2}} \sin \Big(\tfrac{1}{2}\sqrt{1+4a_\ell^2} x_\ell\Big) \)$}\prod_{r=\ell+2}^{2\ell-1} \cos \theta_{r}\Bigg).\end{aligned}\end{equation}

In above, all $a_1,\ldots,a_\ell$ are positive numbers.
\end{theorem}

\section{Type $I$ Hamiltonian-stationary Lagrangian surfaces in $\tilde M^2(4\varepsilon)$.}

For a twisted product decomposition ${\mathcal TP}_{fk}^2( \varepsilon)$ of a simply-connected surface of constant curvature $\varepsilon$,  we have
  \begin{align}  \label{5.1}&\text{\small$\(\frac{f_y}{k}\)_y+\(\frac{k_\ell }{f}\)_x$}=- \varepsilon f k.\end{align} 
The twistor form $\Phi({\mathcal TP})=f^2 dx^2+k^2 dy$ is twisted closed if and only if we have
\begin{align} \label{5.2} ff_y=k k_x.\end{align}
 
 For $n=\ell=2$, Proposition \ref{P:3.1} reduces to

\begin{proposition} \label{P:3.4} \cite{DH} Let $ L_{fk}: {\mathcal TP}_{fk}^2( \varepsilon) \to \tilde M^2(4 \varepsilon)$ be an adapted Lagrangian immersion. Then  $L_{fk}$ is Hamiltonian-stationary if and only if we have  \begin{align}  \label{5.3}& k^3f_x+f^3 k_y=f^2k f_y+fk^2k_x\;\; {\rm (or \; equivalently,} \; \text{\small$  \(\frac{k}{f}\)_x+\(\frac{f}{k}\)_y=0$}).\end{align} 
\end{proposition}

Proposition \ref{P:3.4} implies that the adapted Lagrangian immersion  $L_{fk}$ is always Hamiltonian-stationary   whenever $f^2=k^2$. We call such Hamiltonian-stationary Lagrangian surfaces  to be of {\bf  type $I$.}

It was  proved in \cite{cdvv} that  Hamiltonian-stationary Lagrangian surfaces of  type $I$ in $CP^2(4)$ are  congruent to
\begin{equation}\begin{aligned} &L(x,y)=\text{\small$\frac{1}{\alpha}$} \Big(\text{\small$\frac{ib}{2}$}+\tanh (x+y), e^{{i\over2}b(x+y)}\sech (x+y)\cos(\alpha (x-y)), \\& \hskip.3in e^{{i\over2}b(x+y)}\sech (x+y)\sin(\alpha (x-y))\Big),\; \alpha=\tfrac{1}{2}\sqrt{4+b^2},\; b>0.\end{aligned}\end{equation}  

It is also known in \cite{cdvv} that type $I$ Hamiltonian-stationary Lagrangian immersions in ${\bf C}^2$ are  congruent  to one of
 the following two immersions:
\begin{align}& L=a( e^{ix}, e^{iy}),\;\;\; a>0;
\\&\label{3.17} L=\text{\small$\frac{\sqrt{2}ae^{\frac{1}{2}\sqrt{1+4b^2}(x+y)}}{\sqrt{1+4b^2}} \(\cos 
\(\frac{\sqrt{1+4b^2}}{2}(x-y)\),\sin\(\frac{\sqrt{1+4b^2}}{2}(x-y)\)\)$}.\end{align}

For  type $I$ Hamiltonian-stationary Lagrangian immersions in $CH^2(-4)$,  it is known in \cite{cdvv} that they are congruent to a Lagrangian surfaces obtained from one of the following five families:
\begin{align}&
L=\text{\small$\frac{1}{\alpha}$} \(\text{\small$\frac{ib}{2}$} -\tan s, \text{\small$\frac{e^{{i\over2}bs}\cos(\alpha t)}{\cos s}, \frac{e^{{i\over2}bs}\sin(\alpha t)}{\cos s}$}\),\;
 \alpha^2= {\tfrac{b^2}4} -1, \; b>2; 
 \\ &
L=\text{\small$ \frac{1}{\alpha}$} \(\text{\small$\frac{e^{{i\over2}bs}\cosh(\alpha t)}{\cos s}$} ,\text{\small$\frac{ib}{2}$} -\tan s ,\text{\small$\frac{ e^{{i \over2}bs}\sinh(\alpha t)}{\cos s}$}\),
\,  \alpha^2= 1-{\tfrac{b^2}4}, \, b\in (0,2); 
\\&
L=e^{is}\sec s \(\text{\small$\frac{t^2}{2}+{\frac34}  +\frac{e^{-2is}}{4} -{\frac{is}{2}}$},  t , 
 i  \( \text{\small$\frac{t^2}{2} -{\frac14}  +\frac{e^{-2is}}{4} -{\frac{is}{2}}$}\) \) ;
\\&
L=\text{\small${1\over\alpha}$}\(\text{\small$\frac{ib}{2} $} +\coth s,\text{\small $\frac{e^{{i\over2}bs}\cos(\alpha t)}{\cosh s},\frac{e^{{i\over2}bs}\sin(\alpha t)}{\cosh s}$}\),
\; \alpha^2= {\tfrac{b^2}4} +1, \; b>0; 
\\&
L=\(\text{\small${2\over {x+y}}$}+ i, \text{\small${{\sqrt{2}e^{ix}} \over {x+y}}$},\text{\small${{\sqrt{2}e^{iy}}\over {x+y}}$}\),
\end{align} 
where  $s=x+y$ and $t=x-y$.
\vskip.1in

\section{Type $II$ Hamiltonian-stationary Lagrangian surfaces in $\tilde M^2(4\varepsilon)$.}

In this and next two sections, I  present  my recent joint work with O. J. Garay and Z. Zhou \cite{cgz} concerning   Hamiltonian-stationary Lagrangian surfaces of type {$II$}.
 
If the two twistor functions $f$ and $k$ are unequal, then the twisted product decomposition ${\mathcal TP}_{fk}^2( \varepsilon)$ of a surface of constant curvature $\varepsilon$ gives rise to a $H$-stationary Lagrangian surface in $\tilde M^2(4\varepsilon)$ if and only if $f$ and $k$ satisfy  the  following {\it  over-determined PDE system}: 
  \begin{align}  \label{6.1}& \(\frac{k}{f}\)_x+\(\frac{f}{k}\)_y=0, \\& \label{6.2} \;\;  \frac{k}{f}=\frac{f}{k}, \\    &\label{6.3}\(\frac{f_y}{k}\)_y+\(\frac{k_x}{f}\)_x=-\varepsilon fk,\;\; \hbox{$ (\varepsilon=1, 0$ or $-1$)}.
\end{align}
Hamiltonian-stationary  Lagrangian immersions obtained from such functions $f,k$ with $f^2\ne k^2$ are said to be {\bf  of type II}.

In order to find non-trivial solutions of this over-determined PDE system \e{6.1}-\e{6.3}, first we give the following two lemmas from \cite{cgz}.

\begin{lemma}\label{L:6.1} If $f(x,y)$ and $k(x,y)$ satisfy  Eqs. \e{6.1} and \e{6.2}, then  \begin{align}  \notag F(x,y)=c m f(m^2 x, y),\;\;\; K(x,y)= c k(m^2x, y),\;\; 1\ne m\in {\bf R}^+, \end{align} satisfy Eqs. \e{6.1} and \e{6.2} automatically for any constant $ c \ne 0$.
\end{lemma}
 
The {easiest way to obtain a solution of  system \e{6.1}-\e{6.3} with $f^2\ne k^2$ is to start from  $f,k$ with} $ f^2=k^2$ since \e{6.1} is apparently satisfied, and \e{6.2} and \e{6.3} are in much simpler form. 
  
\begin{lemma}\label{L:6.2} If $f(x,y)$ is a solution of  system \e{6.1}-\e{6.3} with $f^2=k^2$, then  \begin{align}  \label{4.6}f(x,y)= \hat{f} (x+y),  \end{align} for some function $\hat f(x+y)$ which is a traveling wave solution with unit traveling speed. In addition, we have a family of solutions of the system \e{6.1}-\e{6.3} given by
\begin{equation}\begin{aligned} &F(x, y) = \frac{m \sqrt{1+m^2}}{\sqrt{2}} \hat{f}(m^2 x + y), \;\;
\\&  K(x, y) = \pm \frac{ \sqrt{1+m^2}}{\sqrt{2}} \hat{f}(m^2 x + y), \;\; \end{aligned} \end{equation}
where $ m $ is any positive constant.
\end{lemma}

Now, we can apply Lemma \ref{L:6.2} to construct some nontrivial traveling wave solutions of the over-determined PDE system \e{6.1}-\e{6.3} with $f^2 \ne k^2$.

\begin{example} \label{E:6.1}{ \rm When $\varepsilon=1$ and $ f = k = \hat{f}(x+y)$,  Eq. \e{6.3} becomes
\begin{align}\label{6.8} 2 (\ln \hat f)''+ \hat f^2=0.\end{align} 
Since $\hat f \ne 0$,  \e{6.8} implies $\hat f$ is non-constant. Thus,  Eq. \e{6.8}  yields
\begin{align}\label{6.9} 2 \text{\small$\frac{\hat f'{}^2}{\hat f^2}$}+ \hat f^2 =c_1^2 ={\rm constant},\;\; c_1 >0.\end{align} 
After solving \e{6.9}, we know that, up to translations and sign, $\hat f$ is 
\begin{align}\label{6.10}\hat f(u)=c_1 \sech\(\text{\small$\frac{c_1 u}{\sqrt{2}}$} \).\end{align} 
For notational simplicity, taking $ c_1 = \frac{\sqrt{2}c}{\sqrt{1+m^2}} $ and using Lemma \ref{L:6.2}, we obtain the following solutions of system \e{6.1}-\e{6.3}:
\begin{align}\label{6.11}f= c m\sech\(\text{\small$\frac{c(m^2x+y)}{\sqrt{1+m^2}}$} \),\;\; k=\pm c \sech\(\text{\small$\frac{c(m^2x+y)}{\sqrt{1+m^2}}$} \) ,\end{align} 
where $c,m$ are positive real numbers.
}\end{example}

\begin{example}  \label{E:6.2}{\rm When $\varepsilon=0$ and $ f = \pm k = \hat{f}(x+y)$,  Eq. \e{6.3} becomes $(\ln \hat f)''=0$, which implies that  $ \hat f=a e^{bu}$
for some real numbers $a,b$ with $a\ne 0$.
Thus, we have the following solutions of system \e{6.1}-\e{6.3}:
\begin{align}\label{6.13}f(x,y)=a me^{b(m^2x+y)},\;\; k(x,y)=\pm a e^{b(m^2x+y)}, \;\; a\ne 0.\end{align} 
}\end{example}

\begin{example}  \label{E:6.3}{ \rm When $\varepsilon=-1$ and $ f = \pm k = \hat{f} (x+y) $,  \e{6.3} becomes
\begin{align}\label{6.14} 2\(\text{\small$\frac{\hat f'}{\hat f}$}\)'= \hat f^2.\end{align} 
Since $\hat f\ne 0$,  \e{6.14} implies that $\hat f$ is non-constant. Thus, \e{6.14}  yields
\begin{align}\label{6.15} 2 \text{\small$\frac{\hat f'{}^2}{\hat f^2}$}- \hat f^2= \pm c^2 ={\rm constant},\;\; c\geq 0.\end{align} 

{\bf Case} (i): $ - c^2<0 $.  Solving \e{6.15} gives 
$\hat f(u)=c \sec\big(\text{$\frac{cu}{\sqrt{2}}$} \big).$
As in Example \ref{E:6.1}, we find the following solutions of system \e{6.1}-\e{6.3}:
\begin{align}\label{6.17}f(x,y)=c m\sec\Big(\text{\small$\frac{c(m^2x+y)}{\sqrt{1+m^2}}$} \Big),\;\; k(x,y)=\pm c \sec\Big(\text{\small$\frac{c(m^2x+y)}{\sqrt{1+m^2}}$} \Big) ,\end{align} 
where $c,m$ are positive real numbers.
\smallskip

{\bf Case} (ii): $ - c^2<0 $. $ c^2>0 $.  In this case, after solving \e{6.15} we get another family of solutions:
\begin{align} \label{6.18} f(x,y)=c m  \csch \Big(\text{\small$\frac{c(m^2x+y)}{\sqrt{1+m^2}}$} \Big),\;\; k(x,y)=\pm c \csch \Big(\text{\small$\frac{c(m^2x+y)}{\sqrt{1+m^2}}$} \Big) .\end{align} 
\smallskip

{\bf Case} (iii):   $c> 0$.  In this case,  after solving \e{6.15} we get another family of solutions:  
\begin{align} \label{6.19} f(x,y)=\text{\small$ \frac{m\sqrt{1+m^2} }{m^2x+y}$},\;\; k(x,y)=\text{\small$ \frac{\sqrt{1+m^2} }{m^2x+y}$} ,\;\; m\ne 1.\end{align} 
}\end{example}

\section{$H$-stationary surfaces  arisen from traveling wave solutions.}
 We are able to construct  type $II$ Hamiltonian-stationary Lagrangian surfaces in $\tilde M^2(4\varepsilon)$ using  the traveling wave solutions.

First, let us consider the twisted product decomposition  with
\begin{align}\label{7.1} g= \sech^2\(\text{\small$\frac{m^2x+y}{\sqrt{1+m^2}}$} \)(m^2 dx^2+dy^2), \;\;  1\ne m\in {\bf R}.\end{align} 
The horizontal lift $\hat L$ of the
corresponding $H$-stationary Lagrangian immersion $L$
 satisfies the following PDE system:
\begin{equation}\begin{aligned}\label{7.3}&\hat L_{xx}=i \hat L_x- \text{\small$\frac{m^2 (\hat L_x-  \hat L_y) }{\sqrt{1+m^2}} $}\tanh
\( \text{\small$\frac{m^2x+y}{\sqrt{1+m^2}}$} \)
 -m^2  \sech^2\(\hskip-.02in \text{\small$\frac{m^2x+y}{\sqrt{1+m^2}}$}\hskip-.02in \)\hat L,
\\&\hat  L_{xy}=-\frac{\hat L_x+m^2\hat L_y}{\sqrt{1+m^2}}\tanh\(\text{\small$\frac{m^2x+y}{\sqrt{1+m^2}}$} \),\\&
\hat  L_{yy}=i \hat L_y+\text{\small$\frac{ (\hat L_x- \hat L_y)}{\sqrt{1+m^2}}$} \tanh\(\text{\small$\frac{m^2x+y}{\sqrt{1+m^2}}$}\)
- \sech^2\(\text{\small$\frac{m^2x+y}{\sqrt{1+m^2}}$} \)\hat L.\end{aligned}\end{equation} 
Solving this system gives
the following family of   Hamiltonian-stationary Lagrangian surfaces  in $CP^2(4)$:
\begin{equation}\begin{aligned}\notag  &\hat  L=\frac{\sech\, \( \frac{m^2x+y}{\sqrt{1+m^2}} \)}{\text{ $ \sqrt{2+m^2}$}}\Bigg(\text{$\frac{2 m \sqrt{2+m^2}}{\sqrt{1+5m^2}}$}e^{\frac{i }{2}(x+y)} \sin \Big( \frac{\sqrt{1+5m^2}}{2\sqrt{1+m^2}}(x-y)\Big),\\&\hskip0in e^{\frac{i }{2}(x+y)}\hskip-.02in  \left[\text{$ \sqrt{1+m^2}$} \cos \Big( \frac{\sqrt{1+5m^2}}{2\sqrt{1+m^2}}(x-y)\hskip-.02in \Big)\hskip-.02in  - \hskip-.02in \text{$ \frac{i(1-m^2)}{\sqrt{1+5m^2}}$}\sin \Big( \frac{\sqrt{1+5m^2}}{2\sqrt{1+m^2}}(x-y)\hskip-.02in  \Big)\right], \\&\hskip.6in \text{$\frac{1}{\sqrt{2}}\sqrt{1+\cosh \Big(\tfrac{2m^2x+2y}{\sqrt{1+m^2}} \Big)}$}\(1-i\text{$ \sqrt{1+m^2}$}\tanh  \Big(\frac{m^2x+y}{\sqrt{1+m^2}} \Big)\)
\Bigg).
\end{aligned}\end{equation}

Next, let us consider the following  twisted product  flat metric:
\begin{equation}g= e^{2b(m^2x+y)}(m^2 dx^2+dy^2),\; m\ne 0,1; \, b\in {\bf R}.\end{equation} 
The corresponding $H$-stationary Lagrangian immersion $L$
 satisfies the following PDE system:
  \begin{equation}\begin{aligned}&L_{xx}=(i +bm^2)L_x- bm^2  L_y,\; \\&L_{xy}=b L_x+b m^2 L_y,\;
 \\& L_{yy}=-bL_x+(i+b)L_y.\end{aligned}\end{equation} 
Solving this system yields the following new family of flat  $H$-stationary Lagrangian surfaces in {\bf C}$^2$:
\begin{equation}\begin{aligned}&\hskip.0in  L(x,y)=\text{$\frac{e^{\frac{i}{2}(x+y)+b(m^2x+y)}}{\sqrt{1+m^2}}
\Bigg(\frac{2m \sin \big(\frac{1}{2}\sqrt{1+4b^2m^2}(x - y)\big)}{\sqrt{1+4b^2 m^2}},$}
\\&\hskip.2in 
\text{ $ \frac{(1+m^2) \cos \big(\frac{1}{2}\sqrt{1+4b^2m^2}(x - y)\big)}{\sqrt{1+b^2(1+ m^2)^2}}$} \\&\hskip.3in  \text{\small$-\frac{i(1-m^2)\sin \big(\frac{1}{2}\sqrt{1+4b^2m^2}(x - y)\big)}{\sqrt{1+4b^2m^2}\sqrt{1+b^2(1+m^2)^2}}
\Bigg)$}, \end{aligned}\end{equation} where $b,m$ are real numbers with $0<m\ne 1$.
In particular, if $b=0$, we obtain the following new family of $H$-stationary surfaces:
\begin{align}&  L=\text{$\frac{e^{\frac{i}{2}(x+y)}}{\sqrt{1+m^2}}  \sin \(\frac{x - y}{2}\)\Big(2m,
(1+m^2) \cot  \(\frac{x - y}{2}\)-i(1-m^2)\Big)$}. \end{align}

For $CH^2(-4)$ we obtain the following:

\begin{theorem} \cite{cgz} There exist five families of  Hamiltonian-stationary surfaces of type II  in $CH^2(-4)$ arisen the traveling wave solutions:
\begin{equation}\begin{aligned} & {\rm (a)}\; 
\hat L= \Bigg(1-\text{\small$\frac{ i(1+m^2)}{m^2x+y}$} ,\,  \text{\small$\frac{ m\sqrt{1+m^2}}{m^2x+y}$}  e^{ix} ,  \text{\small$\frac{ \sqrt{1+m^2}}{m^2x+y}$} e^{iy}\Bigg); \\&{\rm (b)}\;
\hat L=\text{\small$\sech\, \Big(\frac{x+3y}{2\sqrt{3}} \Big)\( \frac{x-y+4i}{2}e^{\frac{i }{2}(x+y)} , \frac{x-y}{2}e^{\frac{i }{2}(x+y)},\sqrt{3}+2i \tan \Big(\frac{x+3y}{2\sqrt{3}} \Big) \)$};
\\ \notag   &{\rm (c)}\; \hat L=\hskip-.03in \Bigg(\hskip-.03in \frac{\small \text{$ \sqrt{3m^4+2m^2-1}$}\cosh ( \alpha (x-y))+\text{\small$  i(m^2-1)$}\sinh( \alpha (x-y))}{\text{\small$  m \sqrt{3m^2-1} $}e^{-\frac{i }{2}(x+y)}} \sec \big(\frac{m^2x+y}{\sqrt{1+m^2}} \big),
\\& \hskip.3in 
\frac{2m e^{\frac{i }{2}(x+y)} }{\text{\small$  \sqrt{3m^2-1} $}}\sec \Big(\frac{m^2x+y}{\sqrt{1+m^2}} \Big) \sinh ( \alpha (x-y)), \text{$ \frac{1}{m}+ \frac{i\sqrt{1+m^2}}{m}$}\tan\Big(\tfrac{m^2x+y}{\sqrt{1+m^2}} \Big)   \Bigg) ;
 \\  &{\rm (d)}\; \hat L=\hskip-.03in \Bigg(\hskip-.03in \frac{ \text{\small$ \sqrt{1-2m^2-3m^4}$}\cos (\beta (x-y))+\text{\small$  i(1-m^2)$}\sin(\beta (x-y))}{\text{\small$  m \sqrt{1-3m^2} $}e^{-\frac{i }{2}(x+y)} }\sec \Big(\frac{m^2x+y}{\sqrt{1+m^2}} \Big),
\\& \hskip.2in 
\frac{2m e^{\frac{i }{2}(x+y)}}{\text{\small$  \sqrt{1-3m^2} $}} \sec \big(\frac{m^2x+y}{\sqrt{1+m^2}} \big) \sinh (\beta (x-y)), \text{$ \frac{1}{m}+ \frac{i\sqrt{1+m^2}}{m}$}\tan\Big(\frac{m^2x+y}{\sqrt{1+m^2}} \Big)   \Bigg); 
\\ \notag   &{\rm (e)}\; \hat L= \tfrac{1}{\sqrt{2+m^2}}\csch \(\frac{m^2x+y}{\sqrt{1+m^2}}\)\hskip-.03in \Bigg(\hskip-.03in\sinh\(\frac{m^2x+y}{\sqrt{1+m^2}} \)-i\sqrt{1+m^2}\cosh\(\tfrac{m^2x+y}{\sqrt{1+m^2}} \),
\end{aligned}\end{equation} \begin{equation}\begin{aligned}\notag & \hskip.0in e^{\frac{i}{2}(x+y)}\Big\{ \text{$\sqrt{1+m^2}$} \cos \Big(\frac{\sqrt{1+5m^2}}{2\sqrt{1+m^2}}(x-y)\Big)+\frac{i(m^2-1)}{\sqrt{1+5m^2}}\sin \Big(\tfrac{\sqrt{1+5m^2}}{2\sqrt{1+m^2}}(x-y)\Big)\Big\},
\\&\hskip,5in    \frac{2m\sqrt{2+m^2}}{\sqrt{1+5m^2}} e^{\frac{i}{2}(x+y)} \sin \Big(\frac{\sqrt{1+5m^2}}{2\sqrt{1+m^2}}(x-y)\Big)\Bigg), \end{aligned}\end{equation} where $\alpha =\tfrac{\sqrt{3m^2-1}}{2\sqrt{1+m^2}}$ and $\beta =\tfrac{\sqrt{1-3m^2}}{2\sqrt{1+m^2}}$.\end{theorem}

\section{Complete solutions of system \e{6.1}-\e{6.3} for $\varepsilon=0$.}

It is quite difficult to find all exact solutions of the over-determined system. Fortunately,  
we are able to  solve it  for the case $\varepsilon=0$.

\begin{theorem}\label{T:9.1} \cite{cgz} The solutions $\{f,k\}$ of the following over-determined PDE system:
  \begin{align} &\tag{A} \text{ $\(\dfrac{k}{f}\)_x+\(\dfrac{f}{k}\)_y$} =0, 
\;\;\text{ $ \dfrac{f_y}{k}=\dfrac{k_x}{f}, $} \;\; \text{$  \(\dfrac{f_y}{k}\)_y+\(\dfrac{k_x}{f}\)_x$}=0,\end{align}
are  the following:
   \begin{align} \label{8.1}& f(x,y) =\pm k(x,y)=a e^{b(x+y)}; \\ \label{8.2}& f(x,y) =a me^{b(m^2x+y)},\;\;  k(x,y)=\pm a e^{b(m^2x+y)};  \\  \label{8.3}& f(x,y)=\frac{a  \,}{\sqrt{x}}e^{c\arctan \sqrt{-y/x}}, \;\; k(x,y)=\pm \frac{a \, }{\sqrt{-y}}e^{c\arctan \sqrt{-y/x}},  \end{align}
where $a,b,c,m$ are real numbers with $a,c,m\ne 0$ and $m\ne \pm 1$.
\end{theorem}

The  $H$-stationary Lagrangian surfaces of type $II$ in ${\bf C}^2$ corresponding to solutions \e{8.3} 
can also be completely determined as follows: 
 \begin{equation}\begin{aligned}\notag &\hskip.0inL= \text{$ \frac{\sqrt{2\pi} a r^2}{ \sqrt{\cosh(c\pi/2)}}\Bigg($} i J_{-\frac{1}{2}(1+ic)}(r^2)T_c^+(r,\theta)\hskip-.02in +\hskip-.02in  J_{\frac{1}{2}(1-ic)}(r^2)T_c^-(r,\theta) \\&\hskip.4in  +\frac{1}{r^2}\int_0 ^r r e^{ir^2}J_{-\frac{1}{2}(1+ic)}(r^2) dr +\frac{i}{r^2}\int_0^r r e^{ir^2}J_{\frac{1}{2}(1-ic)}(r^2)dr,
\\&\hskip.5in 
ir^2  J_{\frac{1}{2}(1+ic)}(r^2) T_c^+(r,\theta)  
 -r^2  J_{\frac{1}{2}(ic-1)}(r^2)T_c^-(r,\theta)
\\&\hskip.2in +\frac{1}{r^2}\int_0 ^r\hskip-.04in  r e^{ir^2}J_{\frac{1}{2}(1+ic)}(r^2) dr -\frac{i}{r^2}\int_0^r \hskip-.04in   r e^{ir^2}J_{\frac{1}{2}(ic-1)}(r^2)dr\Bigg)
\end{aligned}\end{equation} 
for $a>0,\, c\ne 0$, where  $r,\theta$ and $x,y$ are related by $$x=2r^2 \cos^2\theta, \quad y=-2r^2\sin^2\theta,$$ and $J_\nu (z)$  is the {\it Bessel function of the first kind} with index $\nu$,  which can be expressed the following infinite series:
\begin{align}\label{9.14}J_\nu (z)=\text{\small$ \(\frac{z}{2}\)^\nu \sum_{j=0}^\infty \frac{(-1)^j (z/2)^{2j}}{j! \, \Gamma(\nu+j+1)}.$}\end{align}.

\section{Some applications.}

As an application of the results in section 7, we mention the following classification result from \cite{cg}.

\begin{theorem} \label{T:4.1} There exist five families of Hamiltonian-stationary Lagrangian submanifolds of constant curvature with positive relative nullity in the complex projective $3$-space $CP^3(4)$:
\vskip.04in

{\rm (1)} A totally geodesic Lagrangian submanifold given by $L: RP^3(1)\to CP^3$; 
\vskip.04in

 {\rm (2)}  A Lagrangian submanifold  defined by
\begin{equation}\begin{aligned}\notag &\hskip.4in L(x,y,s)=\Bigg(\text{ $\frac{ \phi e^{\frac{i}{2}s} \sin\(\delta s\) }{ \delta }$}, \text{ $
\frac{\phi e^{\frac{i}{2}s}\{2\delta\cos\(\delta s\)-i \sin\(\delta s\)\}}{\delta\sqrt{4a^2+\hat c^2} }$},\\&\hskip.3in \text{\ $\frac{2\hat c x+ib(1-x^2-y^2)}{\hat c (1+x^2+y^2)}-\frac{2b(\hat c+2ia)\phi}{\hat c(4a^2+\hat c^2)},\,
\frac{2\hat c y+ic(1-x^2-y^2)}{\hat c(1+x^2+y^2)}-\frac{2c(\hat c+2ia)\phi}{\hat c(4a^2+\hat c^2)}$}\Bigg),
\end{aligned}\end{equation}
where $a,b,c$ are real numbers and $$\hat c=\sqrt{b^2+c^2}\ne 0,\;\delta=\frac{1}{2}\sqrt{1+4a^2+\hat c^2},\;\phi=\frac{a(1-x^2-y^2)+bx+cy}{1+x^2+y^2};$$
\vskip.04in

{\rm (3)}   A Lagrangian submanifold defined by
   \begin{equation}\begin{aligned}\notag &\hskip.0in L(x,s,t)=\cos x \(\tan x, \frac{2b\tanh s+i\sqrt{2b}}{\sqrt{2+4b}} ,\frac{\sqrt{2b}e^{{is}/{\sqrt{2b}}} \sech\, s }{\sqrt{1+2b}}
 \cos\(\text{\small$ \frac{\sqrt{1+2b}}{\sqrt{2b}}$}t\)   ,\right. \\&\left. \hskip1.2in \frac{\sqrt{2b}e^{{is}/{\sqrt{2b}}} \sech\, s}{\sqrt{1+2b}} \sin\(\text{\small$ \frac{\sqrt{1+2b}}{\sqrt{2b}}$}t\)  \),\end{aligned}\end{equation}
   where $b$ is a positive number;

{\rm (4)}   A Lagrangian submanifold defined by 
 \begin{equation}   \begin{aligned}  \notag &L(x,y,z)=\Bigg(e^{\frac{iy}{2}}\(\cos \(\tfrac{1}{2}\sqrt{1+4a^2}y\)+\frac{i}{\sqrt{1+4a^2}}\sin\(\tfrac{1}{2}\sqrt{1+4a^2}y\)\)\cos x,\\& \hskip.3in  \frac{2a e^{\frac{iy}{2}}\cos x}{\sqrt{1+4a^2}}\sin\(\tfrac{1}{2}\sqrt{1+4a^2}y\),\frac{2b e^{\frac{iz}{2}}\cos x}{\sqrt{1+4b^2}}\sin\(\tfrac{1}{2}\sqrt{1+4b^2}z\),
   \\& \hskip.4in e^{\frac{iz}{2}}\(\cos \(\tfrac{1}{2}\sqrt{1+4b^2}z\)+\frac{i}{\sqrt{1+4b^2}}\sin\(\tfrac{1}{2}\sqrt{1+4b^2}z\)\)\sin x,\end{aligned}\end{equation}
where $a,b$ are positive numbers.

{\rm (5)} A Lagrangian submanifold defined by
$$L(x,y,z)=\Big(\sin x,\tilde L(y,z)\cos x\Big),$$
where  $\tilde L$ is a horizontal lift of a type II Hamiltonian-stationary Lagrangian surface  $L:{\mathcal TP}_{fk}^2(1)\to CP^2(4)$.
\vskip.04in 

Conversely, locally every  Hamiltonian-stationary Lagrangian submanifold of constant curvature in $CP^3$ with positive relative nullity is   congruent to an open portion of a  Lagrangian submanifold from one of the above five  families.
\end{theorem}

For Hamiltonian-stationary Lagrangian submanifolds in $CH^3$, we have the following result from \cite{c5}.

\begin{theorem} There exist ten families of Hamiltonian-stationary Lagrangian submanifolds of constant curvature in  $CH^3(-4)$  with positive relative nullity: 
\vskip.04in

{\rm (1)} A totally geodesic Lagrangian submanifold $L:H^3(-1)\to CH^3(-4)$; 
\vskip.04in

 {\rm (2)}  A Lagrangian submanifold  defined by
\begin{equation}\begin{aligned}\notag &\hskip.1in L(s,y,z)=\frac{1}{2(1-y^2-z^2)} \Bigg((2i+s)e^{\frac{i}{2}s}\big(2by+\sqrt{1+b^2}(1+y^2+z^2)\big),  \\&\; se^{\frac{i}{2}s}\big(2by+\sqrt{1+b^2}(1+y^2+z^2)\big),4\sqrt{1+b^2}y+2b(1+y^2+z^2), 4z \Bigg),\; b\in {\bf R}.\end{aligned}\end{equation} 
\vskip.04in

{\rm (3)}   A Lagrangian submanifold defined by
  \begin{equation}\begin{aligned}\notag &L(s,y,z)=\frac{1}{1-y^2-z^2}\Bigg(\frac{e^{\frac{i}{2}s}(by+a(1+y^2+z^2))\{2\delta \cosh \delta s-i\sinh \delta s\}}{\delta\sqrt{4a^2-b^2}}\\&\hskip.4in \frac{e^{\frac{i}{2}s}(by+a(1+y^2+z^2)) \sinh \delta s}{\delta}, \frac{4ay-b(1+y^2+z^2)}{\sqrt{4a^2-b^2}},2z\Bigg),\end{aligned}\end{equation} 
   where $a,b,\delta$ are real numbers satisfying $4a^2-b^2>1$ and  $2\delta=\sqrt{4a^2-b^2-1}$.
  \vskip.04in 

{\rm (4)}   A Lagrangian submanifold defined by 
\begin{equation}\begin{aligned}\notag &L(s,y,z)=\frac{1}{1-y^2-z^2}\Bigg(\frac{e^{\frac{i}{2}s}(by+a(1+y^2+z^2))\{2\gamma \cos \gamma s-i\sin \gamma s\}}{\gamma\sqrt{4a^2-b^2}}\\&\hskip.4in \frac{e^{\frac{i}{2}s}(by+a(1+y^2+z^2)) \sin \gamma s}{\gamma}, \frac{4ay+b(1+y^2+z^2)}{\sqrt{4a^2-b^2}},2z\Bigg),\end{aligned}\end{equation} 
  where $a,b,\gamma$ are real numbers satisfying $4a^2<1+b^2$,  $2\gamma=\sqrt{1+b^2-4a^2}$ and $4a^2\ne b^2$.
\vskip.04in 
 
{\rm (5)} A Lagrangian submanifold defined by
\begin{equation}\begin{aligned}\notag &\hskip.0in L(s,y,z)=\Bigg( \frac{2y-a^2(1+is)((1+y)^2+z^2)}{\sqrt{a^2-1}(1-y^2-z^2)},\frac{2z}{1-y^2-z^2},\\& \hskip.2in \frac{1+y^2+z^2+ia^2 s((1+y)^2+z^2)}{\sqrt{a^2-1}(1-y^2-z^2)},\frac{ae^{is}((1+y)^2+z^2)}{1-y^2-z^2}\Bigg),\;\; a^2\ne 0,1. \end{aligned}\end{equation} 
\vskip.04in

{\rm (6)} A Lagrangian submanifold defined by
\begin{equation}\begin{aligned}\notag &\hskip.0in L(s,y,z)=\Bigg(\frac{is}{2}+\frac{3}{2}-i+ \frac{2i-3-is+(2i-2-is)y}{1-y^2-z^2},\frac{2z}{1-y^2-z^2},\\& \hskip.2in \frac{is}{2} -\frac{1}{2}-i+ \frac{1+2i-is+(2+2i-is)y}{1-y^2-z^2},\frac{e^{is}((1+y)^2+z^2)}{1-y^2-z^2}\Bigg). \end{aligned}\end{equation} 
\vskip.04in

{\rm (7)} A Lagrangian submanifold defined by
   \begin{equation}\begin{aligned}\notag &\hskip.0in L(x,s,t)=\frac{\cosh x}{\sqrt{1-2b}}\Bigg(\sqrt{2b}\tan s-i,\sqrt{2b}e^{{is}/{\sqrt{2b}}}\sec s \cos \(\text{\small$\frac{\sqrt{1-2b}}{\sqrt{2b}}$}t\),
   \\&\hskip.1in  \sqrt{2b}e^{{is}/{\sqrt{2b}}}\sec s \sin\(\text{\small$\frac{\sqrt{1-2b}}{\sqrt{2b}}$}t\), \sqrt{1-2b}\tanh x\Bigg),\;\; 0<2b<1.\end{aligned}\end{equation}
\vskip.04in

 {\rm (8)} A Lagrangian submanifold defined by  \begin{equation}\begin{aligned}\notag &\hskip.0in  L(x,s,t)=\frac{\cosh x}{\sqrt{1-2b}}\Bigg(\sqrt{2b}e^{{is}/{\sqrt{2b}}}\sec s \cosh \(\text{\small$ \frac{\sqrt{2b-1}}{\sqrt{2b}}t$}\),\sqrt{2b}\tan s-i,
   \\&\hskip.2in  \sqrt{2b}e^{{is}/{\sqrt{2b}}}\sec s \sinh  \(\text{\small$ \frac{\sqrt{2b-1}}{\sqrt{2b}}t$}\), \sqrt{2b-1}\tanh x\Bigg),\;\; 2b>1.\end{aligned}\end{equation}
\vskip.04in

{\rm (9)} A Lagrangian submanifold defined by
   \begin{equation}\begin{aligned}\notag &\hskip.0in L(x,s,t)=\frac{\cosh x}{\sqrt{2}(1+e^{2is})}\Big(i+2e^{2is}(s+i+it^2), i+2e^{2is}(s+it^2),\\&\hskip1.0in  \sqrt{2}(1+e^{2is})\tanh x,2\sqrt{2} e^{2is}t\Big).\end{aligned}\end{equation}
\vskip.04in

{\rm (10)} A Lagrangian submanifold defined by$$L(x,y,z)=\big(\tilde P(y,z)\cosh x,\sinh x\big),$$
where  $\tilde P$ is a horizontal lift of a type II Hamiltonian-stationary Lagrangian surface  $L:{\mathcal TP}_{f^2 k^2}^n(-1)\to CH^2(-4)$ via the Hopf fibration $\pi:H^{5}_1(-1)\to CH^2(-4)$.

\vskip.04in
Conversely, locally every  Hamiltonian-stationary Lagrangian submanifold of constant curvature in $CH^3(-4)$ with positive relative nullity is   congruent to an open portion of a  Lagrangian submanifold from one of the above tex  families.
\end{theorem}

\begin{remark} ({\bf Added on July 13, 2013}) The PDE system (A) given in Theorem 8.1 was completely solved in \cite{c12}.  In particular, it was proved in \cite{c12} that the PDE system  admits only traveling wave solutions, whenever $\varepsilon\ne 0$.

\end{remark}

\enddocument